\documentclass[review,hidelinks,onefignum,onetabnum]{siamart220329}
\usepackage{amsmath}
\usepackage{lipsum}
\usepackage{amsfonts}
\usepackage{graphicx}
\usepackage{epstopdf}
\usepackage{algorithmic}
\usepackage{amssymb}
\usepackage{amsopn}

\usepackage{lineno}
\nolinenumbers

\ifpdf
\DeclareGraphicsExtensions{.eps,.pdf,.png,.jpg}
\else
\DeclareGraphicsExtensions{.eps}
\fi

\newsiamremark{remark}{Remark}
\newsiamremark{hypothesis}{Hypothesis}
\crefname{hypothesis}{Hypothesis}{Hypotheses}
\newsiamthm{claim}{Claim}
\headers{Convexification for the 3D Problem of Travel Time Tomography}
{M. V. Klibanov, J. Li, V. G. Romanov, and Z. Yang}

\ifpdf
\hypersetup{
  pdftitle={Convexification for the 3D Problem of Travel Time Tomography},
  pdfauthor={M. V. Klibanov, J. Li, V. G. Romanov, and Z. Yang}
}
\fi
\begin{document}

\title{Convexification for the 3D Problem of Travel Time Tomography}
\author{Michael V. Klibanov \thanks{
Department of Mathematics and Statistics, University of North Carolina at
Charlotte, Charlotte, NC, 28223, USA (mklibanv@charlotte.edu)} \and Jingzhi
Li \thanks{
Department of Mathematics \& National Center for Applied Mathematics
Shenzhen \& SUSTech International Center for Mathematics, Southern
University of Science and Technology, Shenzhen 518055, P.~R.~China
(li.jz@sustech.edu.cn)} \and Vladimir G. Romanov \thanks{%
Sobolev Institute of Mathematics, Novosibirsk 630090, Russian Federation
(romanov@math.nsc.ru) } \and Zhipeng Yang \thanks{%
School of Mathematics and Statistics, Lanzhou University, Lanzhou 730000, P.
R. China (yangzp@lzu.edu.cn)} }
\maketitle

\begin{abstract}
The travel time tomography problem is a coefficient inverse problem for the
eikonal equation. This problem has well known applications in seismic. The
eikonal equation is considered here in the circular cylinder, where point
sources run along its axis and measurements of travel times are conductes on
the whole surface of this cylinder. A new version of the globally convergent
convexification numerical method for this problem is developed. Results of
numerical studies are presented.
\end{abstract}

%\title{Convexification for the 3D Problem of Travel Time Tomography
%\thanks{}}



\section{Introduction}

\label{sec:1}

The travel time tomography problem (TTTP) has a long history since it has
broad applications in seismic studies of Earth, see, e.g. \cite{Rom,Vol}.
First, we refer here to the pioneering works of Herglotz in 1905 \cite{H}
and Wiechert and Zoeppritz in 1907 \cite{W}, who have solved this problem in
the 1d case. There are some important uniqueness and stability results for
the TTTP. Since we are focused here on the numerical aspect, we cite only a
few of those \cite{Mukh,MR,PU,Rom}. We also refer to \cite{Sch,Z} for
numerical studies of the TTTP by the methods, which are significantly
different from the convexification method of this paper.

Our goal here is to present the theory and computational experiments for a
new version of the globally convergent convexification numerical method for
the TTTP in the 3d case. The theory of the first version of the
convexification for this problem was published in \cite{kin1}, \cite[Chapter
11]{KL} and numerical studies were published in \cite{kin3d}. However, it
was assumed in these references that the domain of interest $\Omega $ is a
vertically elongated rectangular prism, and point sources run along an
interval of a straight line, which is located in a horizontal plane
\textquotedblleft below" $\Omega .$ As a result, only horizontally oriented
abnormalities were imaged in \cite{kin3d}.

Since we do not want to be restricted to only horizontally oriented
abnormalities, then we consider in this paper the case when $\Omega $ is a
vertically oriented circular cylinder with point sources running along its
axis and measurements of travel times being conducted on the lateral
surface, top and bottom boundaries of this cylinder, see Figure 1. This
source/detector configuration has applications in seismic imaging of wells.
Naturally, our source/detector configuration requires some changes in the
theory of \cite{kin1,kin3d}.

\begin{figure}[tbph]
\centering
\includegraphics[width = 3in]{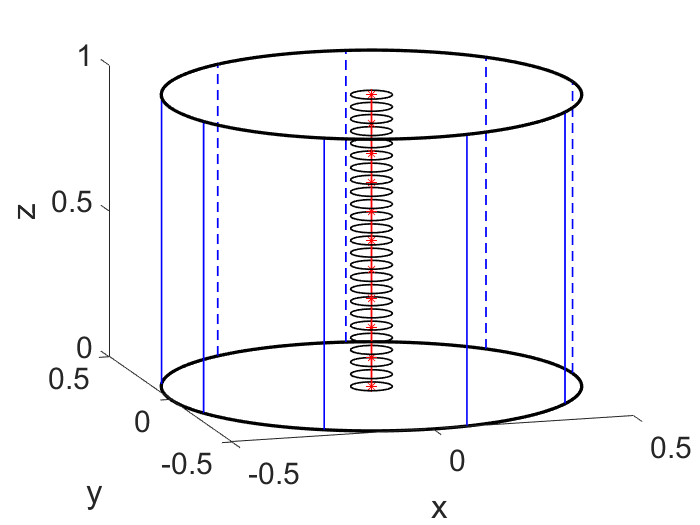}
\caption{A schematic diagram of the source/detectors of this paper. The
large circle is the outer boundary of the domain $\Omega$ defined in (%
\protect\ref{1.3}), the small circle is the inner boundary of domain $%
\Omega_{\protect\varepsilon}$ defined in (\protect\ref{1.31} ), and the red
line with red star is position of the source $L_{\text{source}}$ defined in (%
\protect\ref{1.32}). }
\label{Plot_Domain}
\end{figure}

Below we first formulate and prove a new theorem about the monotonicity of
the solution of the eikonal equation with respect to the radius of the
cylindrical coordinate system. Next, using this theorem, we present a new
version of the convexification method, which is adapted to the
sources/detectors configuration of Figure 1. Next, we formulate some
theorems of the convergence analysis. Since their close analogs were proven
in \cite{kin3d}, then we do not prove them here. Finally, we describe our
numerical studies and present their results.

\textbf{Definition.} \emph{We call a numerical method for a Coefficient
Inverse Problem (CIP) globally convergent, if there exists a theorem
claiming that this method delivers points in a sufficiently small
neighborhood of the true solution of that problem without any advanced
knowledge of this neighborhood. The size of this neighborhood depends only
on the noise level in the input data for this CIP. In other words, a good
first guess about the true solution is not required. }

The first theoretical results about the convexification method for CIPs were
obtained in \cite{KI,Klib97}. More recent works on this method for various
CIPs contain both theory and numerical results, see, e.g. \cite%
{Baud1,Baud2,KLZ,KL,kin3d,convMFG} and references cited therein. The
development of the convexification is caused by the desire to avoid the
phenomenon of local minima and ravines of conventional least squares cost
functionals for CIPs, see, e.g. \cite{BL,Gonch,Grote}.\ In the
convexification, a globally strongly cost functional is constructed. The key
element of this functional is the presence of the Carleman Weight Function
(CWF) in it. This is the function, which is used as the weight in the
Carleman estimate for an associated operator. Usually this is a partial
differential operator. In the case of the TTTP, however, this is a
Volterra-like integral operator. Finally, it is proven that the gradient
descent method of the minimization of this functional converges globally to
the true solution of the corresponding CIP as long as the level of noise in
the data tends to zero.

All functions considered below are real valued ones. In section 2 we pose
our inverse problem. In section 3 we prove the above mentioned important
theorem about the monotonicity of the travel time function with respect to
the radius of the cylindrical coordinate system. In section 4 we derive an
integral differential equation, which does not contain the unknown
coefficient. In section 5 we present an important element of our method: the
orthonormal basis of \cite{KJIIP}, \cite[section 6.2.3]{KL}. In section 6 we
derive a boundary value problem and also prove a Carleman estimate for the
Volterra integral operator. In section 7 we introduce partial finite
differences for our problem. In section 8 we present the convexification
functional and formulate theorems of its global convergence analysis.
Section 9 is devoted to numerical studies.

\section{Statement of the Problem}

\label{sec:2}

Below $\mathbf{x}=\left( x,y,z\right) $ denotes points in $\mathbb{R}^{3}.$
Let $c\left( \mathbf{x}\right) $ be the speed of waves propagation, $c\left( 
\mathbf{x}\right) =1/n\left( \mathbf{x}\right) $, where $n\left( \mathbf{x}%
\right) $ is the refractive index. The function $n\left( \mathbf{x}\right) $
generates the Riemannian metric \cite[Chapter 3]{Rom}%
\begin{equation}
d\tau =n\left( \mathbf{x}\right) \sqrt{\left( dx\right) ^{2}+\left(
dy\right) ^{2}+\left( dz\right) ^{2}}.  \label{1.1}
\end{equation}%
Let $\mathbf{x}_{0}$ be a source of waves and $\mathbf{x}$ be an observation
point. The time, which the wave needs to propagate from a source $\mathbf{x}%
_{0}$ to a point $\mathbf{x,}$ is called the \textquotedblleft first arrival
time" or \textquotedblleft travel time". In the case of a heterogeneous
medium with $n\left( \mathbf{x}\right) \neq const.,$ the first arriving
signal, which arrives at the arrival time to the point $\mathbf{x}$, travels
along the so-called \textquotedblleft geodesic line" $\Gamma \left( \mathbf{x%
},\mathbf{x}_{0}\right) $ generated by metric (\ref{1.1}). The first arrival
time $\tau \left( \mathbf{x},\mathbf{x}_{0}\right) $ is%
\begin{equation*}
\tau \left( \mathbf{x},\mathbf{x}_{0}\right) =\int\limits_{\Gamma \left( 
\mathbf{x},\mathbf{x}_{0}\right) }n\left( \mathbf{y}\left( s\right) \right)
d\sigma ,
\end{equation*}%
where $d\sigma $ is the element of the euclidean length. The function $\tau
\left( \mathbf{x},\mathbf{x}_{0}\right) $ satisfies the following eikonal
equation \cite[Chapter 3]{Rom}:%
\begin{equation}
\begin{split}
& \hspace{1cm}\left( \nabla _{\mathbf{x}}\tau \right) ^{2}=n^{2}\left( 
\mathbf{x}\right) , \\
& \tau \left( \mathbf{x},\mathbf{x}_{0}\right) =O\left( \left\vert \mathbf{x}%
-\mathbf{x}_{0}\right\vert \right) \text{ as }\mathbf{x}\rightarrow \mathbf{x%
}_{0}.
\end{split}
\label{1.2}
\end{equation}

Consider two numbers $R,B>0$ and let $\varepsilon \in \left( 0,R\right) $ be
a small number. We define the domains $\Omega ,\Omega _{\varepsilon }\subset 
\mathbb{R}^{3}$ as:%
\begin{eqnarray}
& \Omega =\left\{ \mathbf{x}=\left( x,y,z\right) :\sqrt{x^{2}+y^{2}}<R,z\in
\left( 0,B\right) \right\} ,  \label{1.3} \\
& \Omega _{\varepsilon } =\left\{ \mathbf{x}=\left( x,y,z\right) :\sqrt{
x^{2}+y^{2}}<\varepsilon ,z\in \left( 0,B\right) \right\} ,  \label{1.31} \\
& \hspace{0.2 cm} L_{\text{source}} =\left\{ \mathbf{x}_{0}=\left(
x_{0},y_{0},z_{0}\right) :x_{0}=y_{0}=0,z_{0}\in \left( 0,B\right) \right\} .
\label{1.32}
\end{eqnarray}%
Hence, $\Omega $ and $\Omega _{\varepsilon }$ are vertically oriented
cylinders with the radii $R$ and $\varepsilon $ respectively, $\Omega
_{\varepsilon }\subset \Omega ,$ and $L_{\text{source}}$ is the vertical
axis of both cylinders, along which sources $\mathbf{x}_{0}$ run.

Consider the cylindrical coordinates $\mathbf{x}=\left( x,y,z\right) =\left(
r\cos \varphi ,r\sin \varphi ,z\right) .$ Therefore by (\ref{1.3}) and (\ref%
{1.31}) in cylindrical coordinates 
\begin{eqnarray}
&&\Omega =\left\{ \left( r,\varphi ,z\right) :r\in \left( 0,R\right)
,\varphi \in \lbrack 0,2\pi ),z\in \left( 0,B\right) \right\} ,  \label{1.4}
\\
&&\Omega _{\varepsilon }=\left\{ \left( r,\varphi ,z\right) :r\in \left(
0,\varepsilon \right) ,\varphi \in \lbrack 0,2\pi ),z\in \left( 0,B\right)
\right\} ,  \label{1.5} \\
&&\hspace{0.5cm}D=\Omega \diagdown \Omega _{\varepsilon }=\left\{ \left(
r,\varphi ,z\right) :r\in \left( \varepsilon ,R\right) ,\varphi \in \lbrack
0,2\pi ),z\in \left( 0,B\right) \right\} .  \label{1.50}
\end{eqnarray}%
We consider below the eikonal equation (\ref{1.2}) only in the domain $%
\Omega .$ We assume that the refractive index $n(\mathbf{x})$ satisfies the
following conditions: 
\begin{align}
n(\mathbf{x})& =1\text{ in }\Omega _{\varepsilon },  \label{1.6} \\
n(\mathbf{x})& \in C^{1}\left( \overline{\Omega }\right) ,  \label{1.7} \\
n_{r}\left( r,\varphi ,z\right) & \geq 0\text{ in }\Omega .  \label{1.8}
\end{align}%
In cylindrical coordinates eikonal equation (\ref{1.2}) has the form 
\begin{equation}
(\tau _{r})^{2}+\frac{1}{r^{2}}(\tau _{\varphi })^{2}+(\tau
_{z})^{2}=n^{2}(r,\varphi ,z).  \label{1.9}
\end{equation}%
By (\ref{1.32}) denote 
\begin{equation}
\tau \left( \mathbf{x,x}_{0}\right) =\tau \left( \mathbf{x,}z_{0}\right) ,%
\text{ }\mathbf{x}_{0}\in L_{\text{source}},\text{ }z_{0}\in \left(
0,B\right) .  \label{1.10}
\end{equation}

\textbf{Travel Time Tomography Problem (TTTP). }\emph{Let conditions (\ref%
{1.3})-(\ref{1.10}) hold. Find the refractive index }$n\left( \mathbf{x}%
\right) $\emph{\ for }$\mathbf{x}\in D,$\emph{\ assuming that the following
functions }$p\left( \varphi ,z,z_{0}\right) $, $p_{0}\left( r,\varphi
,z_{0}\right) ,p_{B}\left( r,\varphi ,z_{0}\right) $\emph{\ are known:}%
\begin{align}
p\left( \varphi ,z,z_{0}\right) & =\tau \left( R,\varphi ,z,z_{0}\right) ,%
\text{ }\forall \varphi \in \left( 0,2\pi \right) ,\forall z,z_{0}\in \left(
0,B\right) ,  \label{1.11} \\
p_{0}\left( r,\varphi ,z_{0}\right) & =\tau \left( r,\varphi ,0,z_{0}\right)
,\forall r\in \left( \varepsilon ,R\right) ,\text{ }\forall \varphi \in
\left( 0,2\pi \right) ,\forall z_{0}\in \left( 0,B\right) ,  \label{1.12} \\
p_{B}\left( r,\varphi ,z_{0}\right) & =\tau \left( r,\varphi ,B,z_{0}\right)
,\forall r\in \left( \varepsilon ,R\right) ,\forall \varphi \in \left(
0,2\pi \right) ,\forall z_{0}\in \left( 0,B\right) .  \label{1.13}
\end{align}

Everywhere below we assume that the following condition is valid:

\medskip \textbf{Assumption of the Regularity of Geodesic Lines}. \emph{For
any pair of points }$\left( \mathbf{x},\mathbf{x}_{0}\right) \in \overline{
\Omega }\times L_{\text{source}}$ \emph{there exists unique geodesic line }$%
\Gamma \left( \mathbf{x},\mathbf{x}_{0}\right) $\emph{\ connecting these two
points.}

\section{The Monotonicity Property of the Function $\protect\tau \left( r, 
\protect\varphi ,z,z_{0}\right) $}

\label{sec:3}

\bigskip

\textbf{Theorem 3.1.} \emph{Assume that conditions (\ref{1.3})-(\ref{1.8})
are in place. Then the following inequality is valid:}%
\begin{equation}
\tau _{r}\left( r,\varphi ,z,z_{0}\right) \geq \frac{\varepsilon }{\sqrt{%
\varepsilon ^{2}+B^{2}}}=c,\text{ }\forall \left( r,\varphi ,z\right) \in 
\overline{D}.  \label{2.1}
\end{equation}

\textbf{Proof.} Denote $\tau _{r}=q_{1}$, $\tau _{\varphi }=q_{2}$, $\tau
_{z}=q_{3}$. Then the eikonal equation (\ref{1.9})\emph{\ }becomes:\emph{\ } 
\begin{equation}
q_{1}^{2}+\frac{1}{r^{2}}q_{2}^{2}+q_{3}^{2}=n^{2}(r,\varphi ,z).
\label{2.2}
\end{equation}%
Differentiating (\ref{2.2}) with respect to $r,\varphi ,z$ and using the
equalities%
\begin{equation*}
\left. 
\begin{array}{c}
(q_{2})_{r}=(q_{1})_{\varphi },\quad
(q_{3})_{r}=(q_{1})_{z},(q_{1})_{\varphi }=(q_{2})_{r},(q_{3})_{\varphi
}=(q_{2})_{z},(q_{2})_{z}=(q_{3})_{\varphi },%
\end{array}%
\right. 
\end{equation*}%
we obtain that along the geodesic line $\Gamma \left( \mathbf{x},\mathbf{x}%
_{0}\right) $ with $\mathbf{x}_{0}=(0,0,z_{0})$%
\begin{equation}
\frac{dr}{ds}=\frac{q_{1}}{n^{2}},\text{ }\frac{d\varphi }{ds}=\frac{q_{2}}{%
r^{2}n^{2}},\text{ }\frac{dz}{ds}=\frac{q_{3}}{n^{2}},  \label{2.30}
\end{equation}%
\begin{equation}
\frac{dq_{1}}{ds}=\frac{n_{r}}{n}+\frac{q_{2}^{2}}{r^{3}n},\text{ }\frac{%
dq_{2}}{ds}=\frac{n_{\varphi }}{n},\text{ }\frac{dq_{3}}{ds}=\frac{n_{z}}{n},%
\text{ }\frac{d\tau }{ds}=1,  \label{2.31}
\end{equation}%
where $ds$ is the element of the Riemannian length of $\Gamma \left( \mathbf{%
x},\mathbf{x}_{0}\right) $.

Equations (\ref{2.30}), (\ref{2.31}) form a coupled system of ordinary
differential equations. We need to figure out initial conditions for this
system. Consider the geodesic line $\Gamma \left( r,\varphi ,z,\mathbf{x}%
_{0}\right) $ that connects the point $\mathbf{x}_{0}=(0,0,z_{0})$ with the
point $(r,\varphi ,z).$ The tangent vector to this line at the point $%
\mathbf{x}_{0}$ is the vector

$\nu _{0}=(\sin \theta _{0}\cos \varphi _{0}$, $\sin \theta _{0}\sin \varphi
_{0}$, $\cos \theta _{0})$. Hence, the initial conditions for this geodesic
line are: 
\begin{align}
r|_{s=0}=0,\quad & \varphi |_{s=0}=\varphi _{0},\quad z|_{s=0}=z_{0},
\label{2.40} \\
(q_{1})_{s=0}=\sin \theta _{0},\quad & (q_{2})_{s=0}=0,\quad
(q_{3})_{s=0}=\cos \theta _{0}.  \label{2.41}
\end{align}%
Thus, solving Cauchy problem (\ref{2.30})-(\ref{2.41}), we obtain equations
of $\Gamma \left( \left( r,\varphi ,z\right) ,\left( 0,0,z_{0}\right)
\right) $ in the parametric form 
\begin{equation*}
r=r(s,z_{0},\theta _{0},\varphi _{0}),\quad \varphi =\varphi (s,z_{0},\theta
_{0},\varphi _{0}),\quad z=z(s,z_{0},\theta _{0},\varphi _{0}).
\end{equation*}%
Due to the regularity of geodesic lines, there exists a one-to-one
correspondence between points $\left( 0,0,z_{0}\right) $ and $(r,\varphi
,z). $ In particular, this means that for any fixed $z_{0}$ there exists a
one-to-one correspondence between $(r,\varphi ,z)$ and parameters $(s,\theta
_{0},\varphi _{0})$.\textbf{\ } We denote this correspondence as: 
\begin{equation*}
s=s(r,\varphi ,z,z_{0})\equiv \tau (r,\varphi ,z,z_{0}),\quad \theta
_{0}=\theta _{0}(r,\varphi ,z,z_{0}),\quad \varphi _{0}=\varphi
_{0}(r,\varphi ,z,z_{0}).
\end{equation*}%
Moreover, we have along the geodesic line $\Gamma \left( \left( r,\varphi
,z\right) ,\left( 0,0,z_{0}\right) \right) :$ 
\begin{equation}
\frac{dq_{1}}{ds}=\frac{d\tau _{r}}{ds}=\frac{n_{r}}{n}+\frac{q_{2}^{2}}{%
r^{3}n^{2}}\geq 0.  \label{2.5}
\end{equation}%
For any point $\left( r,\varphi ,z\right) \in \overline{\Omega }\diagdown
\Omega _{\varepsilon }$ the geodesic line $\Gamma \left( r,\varphi ,z,%
\mathbf{x}_{0}\right) $ is: 
\begin{equation*}
\Gamma \left( r,\varphi ,z,\mathbf{x}_{0}\right) =\Gamma _{\varepsilon
}\left( r,\varphi ,z,\mathbf{x}_{0}\right) \cup \Gamma _{1}\left( r,\varphi
,z,\mathbf{x}_{0}\right) ,
\end{equation*}%
where 
\begin{equation}
\Gamma _{\varepsilon }\left( r,\varphi ,z,\mathbf{x}_{0}\right) =\Gamma
\left( r,\varphi ,z,\mathbf{x}_{0}\right) \cap \Omega _{\varepsilon }\text{
and }\Gamma _{1}\left( r,\varphi ,z,\mathbf{x}_{0}\right) \subset \overline{D%
}.  \label{2.6}
\end{equation}%
By (\ref{1.6}) $\Gamma _{\varepsilon } \left( r,\varphi ,z,\mathbf{x}%
_{0}\right) $ is a part of a straight line. Let $\left( \varepsilon ,\varphi
^{\ast },z^{\ast }\right) $ be the intersection of $\Gamma _{\varepsilon }(%
\mathbf{x,x}_{0})$ with the lateral boundary $\left\{ r=\varepsilon ,z\in
\left( 0,B\right) \right\} $ of the cylinder $\Omega _{\varepsilon }.$ Then $%
\tau \left( r,\varphi ^{\ast },z^{\ast },z_{0}\right) =\sqrt{r^{2}+\left(
z^{\ast }-z_{0}\right) ^{2}}$ for $r\in \left( 0,\varepsilon \right) .$
Hence,%
\begin{equation}
\tau _{r}\left( \varepsilon ,\varphi ^{\ast },z^{\ast },z_{0}\right) =\frac{%
\varepsilon }{\sqrt{\varepsilon ^{2}+\left( z^{\ast }-z_{0}\right) ^{2}}}%
\geq \frac{\varepsilon }{\sqrt{\varepsilon ^{2}+B^{2}}}.  \label{2.7}
\end{equation}%
Hence, by (\ref{1.8}), (\ref{2.31}) and (\ref{2.7}) for $\left( r,\varphi
,z\right) \in \overline{\Omega }\diagdown \Omega _{\varepsilon }$ 
\begin{align*}
\tau _{r}\left( r,\varphi ,z,z_{0}\right) =\int\limits_{\sqrt{\varepsilon
^{2}+\left( z^{\ast }-z_{0}\right) ^{2}}}^{s}& \left( \frac{n_{r}}{n}+\frac{%
q_{2}^{2}}{r^{3}n^{2}}\right) ds+\frac{\varepsilon }{\sqrt{\varepsilon
^{2}+\left( z^{\ast }-z_{0}\right) ^{2}}}\geq \\
& \geq \frac{\varepsilon }{\sqrt{\varepsilon ^{2}+B^{2}}}=c.\text{ \ }\square
\end{align*}

\section{Integral Differential Equation}

\label{sec:4}

Denote 
\begin{equation}
\left. 
\begin{array}{c}
u\left( r,\varphi ,z,z_{0}\right) =\tau _{r}^{2}\left( r,\varphi
,z,z_{0}\right) ,\ \left( r,\varphi ,z\right) \in D,z_{0}\in \left(
0,B\right) .%
\end{array}%
\right.  \label{3.1}
\end{equation}%
Since by Theorem 3.1 $\tau _{r}\left( r,\varphi ,z,z_{0}\right) \geq c>0,$
then \ by (\ref{3.1})%
\begin{equation}
\tau _{r}\left( r,\varphi ,z,z_{0}\right) =\sqrt{u\left( r,\varphi
,z,z_{0}\right) }.  \label{3.4}
\end{equation}%
Hence, (\ref{1.11}) and (\ref{3.4}) imply 
\begin{equation}
\left. 
\begin{array}{c}
\tau \left( r,\varphi ,z,z_{0}\right) =-\int\limits_{r}^{R}\sqrt{u\left(
t,\varphi ,z,z_{0}\right) }dt+p\left( \varphi ,z,z_{0}\right) , \\ 
\tau \left( r,\varphi ,0,z_{0}\right) =p_{0}\left( r,\varphi ,z_{0}\right)
,\ \tau \left( r,\varphi ,B,z_{0}\right) =p_{B}\left( r,\varphi
,z_{0}\right) , \\ 
\left( r,\varphi ,z,z_{0}\right) \in Q=\left( \varepsilon ,R\right) \times
\left( 0,2\pi \right) \times \left( 0,B\right) \times \left( 0,B\right) .%
\end{array}%
\right.  \label{3.5}
\end{equation}

By the first line of (\ref{3.5})%
\begin{align}
\tau _{\varphi }\left( r,\varphi ,z,z_{0}\right) & =-\frac{1}{2}%
\int\limits_{r}^{R}\frac{u_{\varphi }\left( t,\varphi ,z,z_{0}\right) }{%
\sqrt{u\left( t,\varphi ,z,z_{0}\right) }}dt+p_{\varphi }\left( \varphi
,z,z_{0}\right) ,  \label{3.6} \\
\tau _{z}\left( r,\varphi ,z,z_{0}\right) & =-\frac{1}{2}\int\limits_{r}^{R}%
\frac{u_{z}\left( t,\varphi ,z,z_{0}\right) }{\sqrt{u\left( t,\varphi
,z,z_{0}\right) }}dt+p_{z}\left( \varphi ,z,z_{0}\right) .  \label{3.7}
\end{align}

\textbf{\ }Substituting (\ref{3.5})-(\ref{3.7}) in the eikonal equation (\ref%
{1.9}), we obtain the following equation in the domain $Q$ defined in the
last line of (\ref{3.5}):%
\begin{equation*}
u\left( r,\varphi ,z,z_{0}\right) +\frac{1}{r^{2}}\left[ -\frac{1}{2}%
\int\limits_{r}^{R}\frac{u_{\varphi }\left( t,\varphi ,z,z_{0}\right) }{%
\sqrt{u\left( t,\varphi ,z,z_{0}\right) }}dt+p_{\varphi }\left( \varphi
,z,z_{0}\right) \right] ^{2}+
\end{equation*}%
\begin{equation}
+\left[ -\frac{1}{2}\int\limits_{r}^{R}\frac{u_{z}\left( t,\varphi
,z,z_{0}\right) }{\sqrt{u\left( t,\varphi ,z,z_{0}\right) }}dt+p_{z}\left(
\varphi ,z,z_{0}\right) \right] ^{2}=n^{2}(r,\varphi ,z).  \label{3.8}
\end{equation}
Differentiating (\ref{3.8}) with respect to $z_{0}$ and using $\partial
_{z_{0}}\left[ n^{2}(r,\varphi ,z)\right] \equiv 0,$ we obtain 
\begin{equation*}
u_{z_{0}}\left( r,\varphi ,z,z_{0}\right) +\frac{1}{r^{2}}\frac{\partial }{%
\partial z_{0}}\left[ -\frac{1}{2}\int\limits_{r}^{R}\frac{u_{\varphi
}\left( t,\varphi ,z,z_{0}\right) }{\sqrt{u\left( t,\varphi ,z,z_{0}\right) }%
}dt+p_{\varphi }\left( \varphi ,z,z_{0}\right) \right] ^{2}+
\end{equation*}%
\begin{equation}
+\frac{\partial }{\partial z_{0}}\left[ -\frac{1}{2}\int\limits_{r}^{R}\frac{%
u_{z}\left( t,\varphi ,z,z_{0}\right) }{\sqrt{u\left( t,\varphi
,z,z_{0}\right) }}dt+p_{z}\left( \varphi ,z,z_{0}\right) \right] ^{2}=0.
\label{3.9}
\end{equation}%
In addition, (\ref{1.12}), (\ref{1.13}) and (\ref{3.1}) imply:%
\begin{equation}
u\left( r,\varphi ,0,z_{0}\right) =\left( \partial _{r}p_{0}\left( r,\varphi
,z_{0}\right) \right) ^{2},\text{ }u\left( r,\varphi ,B,z_{0}\right) =\left(
\partial _{r}p_{B}\left( r,\varphi ,z_{0}\right) \right) ^{2}.  \label{3.10}
\end{equation}%
In (\ref{3.8})-(\ref{3.10}) 
\begin{equation}
\left( r,\varphi ,z,z_{0}\right) \in \overline{Q}.  \label{3.11}
\end{equation}

Thus, we came up with the following boundary value problem (BVP).

\textbf{Boundary Value Problem 1 (BVP1).} \emph{Find the function }$u\left(
r,\varphi ,z,z_{0}\right) \in C^{1}\left( \overline{Q}\right) $\emph{\ \
satisfying conditions (\ref{3.9}), (\ref{3.10}) in the domain (\ref{3.11}).}

Suppose that we have solved this problem. Then we substitute its solution in
the left hand side of equation (\ref{3.8}) and find the target function $%
n(r,\varphi ,z).$ Clearly BVP1\ is a very complicated one. Therefore, we
construct below an approximate mathematical model for its solution, and then
confirm the validity of this model computationally.

\section{A Special Orthonormal Basis in $L_{2}\left( 0,B\right) $}

\label{sec:5}

The basis described in this section was first introduced in \cite{KJIIP}.
Then it was used for the convexification method in some other works, see,
e.g. \cite{kin1,KL,kin3d}. Consider the set of functions $%
\{z_{0}^{n}e^{z_{0}}\}_{n=0}^{\infty }\subset L_{2}(0,B).$ Obviously, this
set is complete in $L_{2}(0,B),$ and these functions are linearly
independent. Apply the Gram-Schmidt orthonormalization procedure to this
set. Then we obtain the orthonormal basis $\{\Psi _{n}\left( z_{0}\right)
\}_{n=0}^{\infty }$ in $L_{2}(0,B)$. The function $\Psi _{n}(z_{0})$ has the
form $\Psi _{n}(z_{0})=P_{n}(z_{0})e^{z_{0}}$, $\forall n\geq 0,$ where $%
P_{n}(z_{0})$ is a polynomial of the degree $n$. Denote 
\begin{equation*}
a_{mn}=\int\limits_{0}^{B}\Phi _{m}\left( z_{0}\right) ,\Phi _{n}^{\prime
}\left( z_{0}\right) dz_{0}.
\end{equation*}%
It was proven in \cite{KJIIP}, \cite[Theorem 6.2.1]{KL} that 
\begin{equation}
a_{mn}=\left\{ 
\begin{array}{c}
1\text{ if }m=n, \\ 
0\text{ if }m>n.%
\end{array}%
\right.  \label{4.1}
\end{equation}%
Let $N\geq 1$ be an integer. Consider the $N\times N$ matrix $A_{N}=\left(
a_{mn}\right) _{m,n=0}^{N}.$ By (\ref{4.1}) $\det A_{N}=1.$ Hence, the
matrix $A_{N}$ is invertible. We observe that neither classical orthonormal
polynomials nor the basis of trigonometric functions do not provide a
corresponding invertible matrix $A_{N}$. This is because the first function
in these two cases is an identical constant, implying that the first column
of that analog of $A_{N}$ is formed only by zeros.

\section{Boundary Value Problem 2}

\label{sec:6}

The first step of our approximate mathematical model mentioned in section 4
is the assumption that the function $u$ can be represented as truncated
Fourier-like series, 
\begin{equation}
u\left( r,\varphi ,z,z_{0}\right) =\sum\limits_{s=0}^{N-1}u_{s}\left(
r,\varphi ,z\right) \Phi _{s}\left( z_{0}\right) ,\text{ }\left( r,\varphi
,z,z_{0}\right) \in Q.  \label{5.1}
\end{equation}%
We do not prove convergence of our method as $N\rightarrow \infty .$ We
point out that the absence of such proofs is a rather common place in the
theory of Inverse Problems, basically due to their ill-posedness, see, e.g. 
\cite{GN,Kab}.

Consider the $N-$D vector function $V\left( r,\varphi ,z\right) $ of
coefficients of the series (\ref{5.1}) 
\begin{equation}
\left. 
\begin{array}{c}
V\left( r,\varphi ,z\right) =\left( u_{0},...,u_{N-1}\right) ^{T}\left(
r,\varphi ,z\right) ,\text{ } \\ 
\left( r,\varphi ,z\right) \in \left( \varepsilon ,R\right) \times \left(
0,2\pi \right) \times \left( 0,B\right) .%
\end{array}%
\right.  \label{5.2}
\end{equation}%
We want to find the vector function $V\left( r,\varphi ,z\right) $ in (\ref%
{5.2}). Naturally, we assume that for functions in (\ref{1.11})-(\ref{1.13}) 
\begin{equation}
\left. 
\begin{array}{c}
p\left( \varphi ,z,z_{0}\right) =\sum\limits_{s=0}^{N-1}p_{s}\left( \varphi
,z\right) \Phi _{s}\left( z_{0}\right) ,\quad \\ 
\left( \varphi ,z,z_{0}\right) \in \left( 0,2\pi \right) \times \left(
0,B\right) \times \left( 0,B\right) .%
\end{array}%
\right.  \label{5.3}
\end{equation}%
\begin{equation}
\left. 
\begin{array}{c}
\left( \partial _{r}p_{0}\right) ^{2}\left( r,\varphi ,z_{0}\right)
=\sum\limits_{s=0}^{N-1}\widetilde{p}_{0,s}\left( r,\varphi \right) \Phi
_{s}\left( z_{0}\right) , \\ 
\left( \partial _{r}p_{B}\right) ^{2}\left( r,\varphi ,z_{0}\right)
=\sum\limits_{s=0}^{N-1}\widetilde{p}_{B,s}\left( r,\varphi \right) \Phi
_{s}\left( z_{0}\right) , \\ 
\left( r,\varphi ,z_{0}\right) \in \left( \varepsilon ,R\right) \times
\left( 0,2\pi \right) \times \left( 0,B\right) .%
\end{array}%
\right.  \label{5.4}
\end{equation}%
Denote%
\begin{equation}
\left. 
\begin{array}{c}
G\left( \varphi ,z\right) =\left( p_{0},...,p_{N-1}\right) ^{T}\left(
\varphi ,z\right) ,\quad \left( \varphi ,z\right) \in \left( 0,2\pi \right)
\times \left( 0,B\right) , \\ 
\text{ }G_{0}\left( r,\varphi \right) =\left( \widetilde{p}_{0,0},...%
\widetilde{p}_{0,N-1}\right) ^{T}\left( r,\varphi \right) ,\text{ }\left(
r,\varphi \right) \in \left( \varepsilon ,R\right) \times \left( 0,2\pi
\right) \\ 
G_{B}\left( r,\varphi \right) =\left( \widetilde{p}_{B,0},...\widetilde{p}%
_{B,N-1}\right) ^{T}\left( r,\varphi \right) ,\text{ }\left( r,\varphi
\right) \in \left( \varepsilon ,R\right) \times \left( 0,2\pi \right) .%
\end{array}%
\right.  \label{5.5}
\end{equation}

Substituting (\ref{5.1})-(\ref{5.5}) in (\ref{3.9}), multiplying the
obtained equality sequentially by the functions $\Psi _{n}\left(
z_{0}\right) ,n=0,...,N-1$ and integrating with respect to $z_{0}\in \left(
0,B\right) ,$ we obtain 
\begin{equation}
\left. 
\begin{array}{c}
A_{N}V+S\left( V_{\varphi },V_{z},G,r,\varphi ,z\right) =0, \\ 
V\left( r,\varphi ,0\right) =\text{ }G_{0}\left( r,\varphi \right) , \\ 
V\left( r,\varphi ,B\right) =\text{ }G_{B}\left( r,\varphi \right) , \\ 
V\left( r,0,z\right) =V\left( r,2\pi ,z\right) , \\ 
\left( r,\varphi ,z\right) \in \left( \varepsilon ,R\right) \times \left(
0,2\pi \right) \times \left( 0,B\right) ,%
\end{array}%
\right.  \label{5.7}
\end{equation}%
where $S\in \mathbb{R}^{N}$ is a $C^{1}-$vector function of its $3N+3$
arguments. The fourth line in (\ref{5.7}) means the $2\pi -$periodicity
condition with respect to the angle $\varphi .$ Thus, (\ref{5.7}) is the
boundary value problem for the nonlinear system of integral differential
equations with respect to the vector function $V\left( r,\varphi ,z\right) .$
Here, the $s-$th component of the $N-$D vector function $S=\left(
S_{0},...,S_{N-1}\right) ^{T}$ has the form%
\begin{equation*}
S_{s}\left( V_{\varphi },V_{z},r\right) =S_{s}^{\varphi }\left( V_{\varphi
},V_{z},r\right) +S_{s}^{z}\left( V_{\varphi },V_{z},r\right) .
\end{equation*}%
Here $S_{n}^{\varphi }\left( V_{\varphi },V_{z},G_{\varphi },G_{z},r\right) $
has the following form after the integration by parts with respect to $z_{0}$%
:%
\begin{align}
S_{s}^{\varphi }\left( V_{\varphi },V_{z},r\right) =& \Psi _{s}\left(
B\right) \frac{1}{r^{2}}\left[ -\frac{1}{2}\int\limits_{r}^{R}\frac{%
u_{\varphi }\left( t,\varphi ,z,B\right) }{\sqrt{u\left( t,\varphi
,z,B\right) }}dt+p_{\varphi }\left( \varphi ,z,B\right) \right] ^{2}-  \notag
\\
-\Psi _{s}\left( 0\right) \frac{1}{r^{2}}& \left[ -\frac{1}{2}%
\int\limits_{r}^{R}\frac{u_{\varphi }\left( t,\varphi ,z,0\right) }{\sqrt{%
u\left( t,\varphi ,z,0\right) }}dt+p_{\varphi }\left( \varphi ,z,0\right) %
\right] ^{2}-  \label{5.8} \\
-\frac{1}{r^{2}}\int\limits_{0}^{B}\Psi _{s}^{\prime }& \left( z_{0}\right) %
\left[ -\frac{1}{2}\int\limits_{r}^{R}\frac{u_{\varphi }\left( t,\varphi
,z,z_{0}\right) }{\sqrt{u\left( t,\varphi ,z,z_{0}\right) }}dt+p_{\varphi
}\left( \varphi ,z,z_{0}\right) \right] ^{2} d z_0.  \notag
\end{align}%
And $S_{s}^{z}\left( V_{\varphi },V_{z},G_{\varphi },G_{z},r\right) $ has
the form after the integration by parts with respect to $z_{0}$:%
\begin{align}
S_{s}^{z}\left( V_{\varphi },V_{z},r\right) =& \Psi _{s}\left( B\right) 
\left[ -\frac{1}{2}\int\limits_{r}^{R}\frac{u_{z}\left( t,\varphi
,z,B\right) }{\sqrt{u\left( t,\varphi ,z,B\right) }}dt+p_{z}\left( \varphi
,z,B\right) \right] ^{2}-  \label{5.9} \\
-\Psi _{s}\left( 0\right) & \left[ -\frac{1}{2}\int\limits_{r}^{R}\frac{%
u_{z}\left( t,\varphi ,z,0\right) }{\sqrt{u\left( t,\varphi ,z,0\right) }}%
dt+p_{z}\left( \varphi ,z,0\right) \right] ^{2}-  \notag \\
-\int\limits_{0}^{B}\Psi _{s}^{\prime }& \left( z_{0}\right) \left[ -\frac{1%
}{2}\int\limits_{r}^{R}\frac{u_{z}\left( t,\varphi ,z,z_{0}\right) }{\sqrt{%
u\left( t,\varphi ,z,z_{0}\right) }}dt+p_{z}\left( \varphi ,z,z_{0}\right) %
\right] ^{2}dz_{0}.  \notag
\end{align}%
In (\ref{5.8}), (\ref{5.9}) the function $u$ has the form (\ref{5.1}), and
the vector function $V$ has the form (\ref{5.2}). Thus, we have obtained
Boundary Value Problem 2 (BVP2).

\textbf{Boundary Value Problem 2 (BVP2).} \emph{Let }$D$\emph{\ be the
domain defined in (\ref{1.50}). Find the }$N-$\emph{D vector function }$V\in
C^{1}\left( \overline{D}\right) $\emph{\ satisfying conditions (\ref{5.7}),
where the function }$u$\emph{\ is connected with }$V$\emph{\ via (\ref{5.1}%
), (\ref{5.2}).}

We now formulate a Carleman estimate for the Volterra integral operator
occurring in (\ref{5.7})-(\ref{5.9}). First, we introduce the Carleman
Weight Function $\varphi _{\lambda }\left( r\right) $ for this Volterra
operator, where $\lambda >0$ is a parameter. This function is:%
\begin{equation}
\varphi _{\lambda }\left( r\right) =e^{2\lambda r}.  \label{5.140}
\end{equation}

\textbf{Theorem 6.1 }(Carleman estimate for the Volterra integral operator). 
\emph{The following Carleman estimate is valid:}%
\begin{equation}
\int\limits_{\varepsilon }^{R}\left( \int\limits_{r}^{R}f\left( s\right)
ds\right) ^{2}e^{2\lambda r}dr\leq \frac{1}{\lambda ^{2}}\int\limits_{%
\varepsilon }^{R}f^{2}\left( r\right) e^{2\lambda r}dr,\text{ }\forall
\lambda >0,\text{ }\forall f\in L_{2}\left( \varepsilon ,R\right) .
\label{5.14}
\end{equation}

\textbf{Proof.} Integrating by parts and using Cauchy-Schwarz inequality, we
obtain%
\begin{align*}
&\int\limits_{\varepsilon }^{R}\left( \int\limits_{r}^{R}f\left( s\right)
ds\right) ^{2} e^{2\lambda r}dr=-\frac{1}{2\lambda }\left(
\int\limits_{\varepsilon }^{R}f\left( s\right) ds\right) ^{2}e^{2\lambda
\varepsilon }dr+ \\
& \hspace{1 cm} +\frac{1}{\lambda }\int\limits_{\varepsilon }^{R}\left[
e^{\lambda r}\int\limits_{r}^{R}f\left( s\right) ds\right] \left[ f\left(
r\right) e^{\lambda r}\right] dr\leq \\
&\leq \frac{1}{\lambda }\left[ \int\limits_{\varepsilon }^{R}\left(
\int\limits_{r}^{R}f\left( s\right) ds\right) ^{2}e^{2\lambda r}dr\right]
^{1/2}\left[ \int\limits_{\varepsilon }^{R}f^{2}\left( r\right) e^{2\lambda
r}dr\right] ^{1/2}.
\end{align*}%
Thus, we have proven that 
\begin{equation}
\left. 
\begin{array}{c}
\int\limits_{\varepsilon }^{R}\left( \int\limits_{r}^{R}f\left( s\right)
ds\right) ^{2}\varphi _{\lambda }\left( r\right) dr\leq \\ 
\leq \left( 1/\lambda \right) \left[ \int\limits_{\varepsilon }^{R}\left(
\int\limits_{r}^{R}f\left( s\right) ds\right) ^{2}\varphi _{\lambda }\left(
r\right) dr\right] ^{1/2}\left[ \int\limits_{\varepsilon }^{R}f^{2}\left(
r\right) \varphi _{\lambda }\left( r\right) dr\right] ^{1/2}.%
\end{array}
\right.  \label{5.15}
\end{equation}%
Dividing both sides of (\ref{5.15}) by 
\begin{equation*}
\left[ \int\limits_{\varepsilon }^{R}\left( \int\limits_{r}^{R}f\left(
s\right) ds\right) ^{2}\varphi _{\lambda }\left( r\right) dr\right] ^{1/2}
\end{equation*}%
and then squaring both sides of the resulting inequality, we obtain (\ref%
{5.14}). $\square $

A similar estimate was proven in \cite{kin1}. However, estimate (\ref{5.14})
is stronger than the one of \cite{kin1} due to the presence of the
multiplier $1/\lambda ^{2}$ here instead of $1/\lambda $ in \cite{kin1} and $%
1/\lambda ^{2}<<1/\lambda $ for $\lambda >>1.$ Note that we use $\lambda >>1 
$ below.\ 

\section{Partial Finite Differences}

\label{sec:7}

To solve boundary value problem (\ref{5.7}), we use partial finite
differences with respect to the variables $z$ and $\varphi .$ without
imposing an additional assumption. This is the second and the last step of
our approximate mathematical model. \ More precisely, we assume that problem
(\ref{5.7}) in finite differences with respect to $z$ and $\varphi ,$ where
the grid step sizes $h_{z}$ and $h_{\varphi }$ with respect to $z$ and $%
\varphi $ satisfy the following inequality:%
\begin{equation}
h_{z},h_{\varphi }\geq h_{0}>0,  \label{7.1}
\end{equation}%
where the number $h_{0}$ is fixed. Consider finite difference grids with
respect to $z$ and $\varphi ,$%
\begin{align}
& 0=z^{0}<z_{1}<...<z_{k-1}<z_{k}=B,z_{i}-z_{i-1}=h_{z},  \label{7.2} \\
0& =\varphi _{0}<\varphi _{1}<...<\varphi _{n-1}<\varphi _{n}=2\pi ,\varphi
_{i}-\varphi _{i-1}=h_{\varphi }.  \label{7.3}
\end{align}

Denote $h=\left( h_{z},h_{\varphi }\right) .$ The vector function $V\left(
r,\varphi ,z\right) $ is turned now in the semi-discrete function, which is
defined on the grid (\ref{7.2}), (\ref{7.3}) for all $r\in \left[
\varepsilon ,R\right] ,$ i.e. $r$ remains a continuous variable. We denote
this function $V^{h}\left( r\right) $. The semi-discrete form of (\ref{5.2})
is:%
\begin{equation}
\left. 
\begin{array}{c}
V^{h}\left( r\right) =\left\{ \left( u_{0},...,u_{N-1}\right) ^{T}\left(
r,\varphi _{i},z_{j}\right) \right\} _{\left( i,j\right) =\left( 0,0\right)
}^{\left( n,k\right) },\text{ } \\ 
\left( r,\varphi _{i},z_{j}\right) \in \widetilde{Q}=\left( \varepsilon
,R\right) \times \left( 0,2\pi \right) \times \left( 0,B\right) .%
\end{array}%
\right.  \label{7.30}
\end{equation}%
And the semi-discrete form of (\ref{5.1}) is:%
\begin{equation}
u^{h}\left( r,\varphi _{i},z_{j},z_{0}\right)
=\sum\limits_{s=0}^{N-1}u_{s}\left( r,\varphi _{i},z_{j}\right) \Phi
_{s}\left( z_{0}\right) ,\text{ }\left( r,\varphi _{i},z_{j},z_{0}\right)
\in Q.  \label{7.31}
\end{equation}

\subsection{Some specifics of partial finite differences}

\label{sec:7.1}

We use the second order approximation accuracy for derivatives $\partial
_{z}V^{h},\partial _{\varphi }V^{h},$%
\begin{align}
\partial _{z}V^{h}\left( r,\varphi _{i},z_{j}\right) & =\frac{V^{h}\left(
r,\varphi _{i},z_{j+1}\right) -V^{h}\left( r,\varphi _{i},z_{j-1}\right) }{%
2h_{z}},\text{ }i\in \left[ 0,n\right] ,j\in \left[ 1,k-1\right] ,
\label{7.4} \\
\partial _{\varphi }V^{h}\left( r,\varphi _{i},z_{j}\right) & =\frac{%
V^{h}\left( r,\varphi _{i+1},z_{j}\right) -V^{h}\left( r,\varphi
_{i-1},z_{j}\right) }{2h_{\varphi }},\text{ }i\in \left[ 1,n-1\right] ,j\in %
\left[ 0,k\right] .  \label{7.5}
\end{align}%
A separate question is about ensuring boundary conditions in the second and
third lines of (\ref{5.7}) as well as the periodicity condition in the
fourth line of (\ref{5.7}). To satisfy boundary conditions in the second and
third lines of (\ref{5.7}), we use in (\ref{7.4}): 
\begin{align}
V_{z}^{h}\left( r,\varphi _{i},z_{1}\right) & =\frac{V^{h}\left( r,\varphi
_{i},z_{2}\right) -\text{ }G_{0}\left( r,\varphi _{i}\right) }{2h_{z}},\text{
}i\in \left[ 0,n\right] ,  \label{7.6} \\
V_{z}^{h}\left( r,\varphi _{i},z_{k-1}\right) & =\frac{\text{ }G_{B}\left(
r,\varphi _{i}\right) -\text{ }V^{h}\left( r,\varphi _{i},z_{k-2}\right) }{%
2h_{z}},\text{ }i\in \left[ 0,n\right] .  \label{7.7}
\end{align}%
In addition, to ensure the periodicity condition in the fourth line of (\ref%
{5.7}), we set%
\begin{equation}
V\left( r,\varphi _{0},z_{j}\right) =V\left( r,\varphi _{n},z_{j}\right) ,%
\text{ }j\in \left[ 0,k\right]  \label{7.8}
\end{equation}%
with corresponding changes in formulas (\ref{7.4})-(\ref{7.7}) for $i=0,n.$

Let $H$ be a Hilbert space. Consider the direct product of $N$ such spaces $%
H_{N}=H\times H\times ...H,$ $N$ times. If $\left\Vert U\right\Vert _{H}$ is
the norm of $U\in H$ in $H$, then the norm of the vector function $%
\widetilde{U}=\left( U_{0},...,U_{N-1}\right) ^{T}\in H_{N}$ is set below as:%
\begin{equation*}
\left\Vert \widetilde{U}\right\Vert _{H_{N}}=\left[ \sum\limits_{i=0}^{N-1}%
\left\Vert U_{i}\right\Vert _{H}^{2}\right] ^{1/2}.
\end{equation*}%
We now introduce some Hilbert spaces of semi-discrete functions, similarly
with \cite{kin3d}. Denote $\mathbf{x}^{h}=\left\{ \left( r,\varphi
_{i},z_{j}\right) _{\left( i,j\right) =\left( 0,0\right) }^{\left(
n,k\right) },r\in \left( \varepsilon ,R\right) \right\} .$ Along with the
above notation $V^{h}(r$, $\varphi _{i}$, $z_{j}),$ we also use the
equivalent one, which is $V^{h}\left( \mathbf{x}^{h}\right) .$ For the
domain $D$ defined in (\ref{1.50}) we denote%
\begin{equation*}
D^{h}=\left\{ \left( r,\varphi _{i},z_{j}\right) _{\left( i,j\right) =\left(
0,0\right) }^{\left( n,k\right) },\text{ }r\in \left( \varepsilon ,R\right)
\right\} .
\end{equation*}%
Next, using (\ref{7.2}) and (\ref{7.3}), we define those Hilbert spaces as:%
\begin{equation}
\left. 
\begin{array}{c}
L_{2,N}^{h}\left( D^{h}\right) =\left\{ 
\begin{array}{c}
V^{h}(\mathbf{x}^{h}):\left\Vert V^{h}(\mathbf{x}^{h})\right\Vert
_{L_{2,N}^{h}\left( D^{h}\right) }^{2}= \\ 
=\sum\limits_{\left( i,j\right) =\left( 0,0\right) }^{\left( n,k\right)
}\left\{ \int\limits_{\varepsilon }^{R}\left( V^{h}\left( r,\varphi
_{i},z_{j}\right) \right) ^{2}dr\right\} <\infty%
\end{array}%
\right\} , \\ 
H_{N}^{1,h}\left( D^{h}\right) =\left\{ 
\begin{array}{c}
V^{h}(\mathbf{x}^{h}):\left\Vert V^{h}\right\Vert _{H_{N}^{1,h}\left(
D^{h}\right) }^{2}= \\ 
=\left\Vert V_{z}^{h}\right\Vert _{L_{2N}^{h}\left( D^{h}\right)
}^{2}+\left\Vert V_{\varphi }^{h}\right\Vert _{L_{2N}^{h}\left( D^{h}\right)
}^{2}+\left\Vert V^{h}\right\Vert _{L_{2N}^{h}\left( D^{h}\right)
}^{2}<\infty%
\end{array}%
\right\} ,%
\end{array}%
\right.  \label{7.9}
\end{equation}%
\begin{equation}
H_{0,N}^{1,h}\left( D^{h}\right) =\left\{ 
\begin{array}{c}
V^{h}(\mathbf{x}^{h})\in H_{N}^{1,h}\left( D^{h}\right) : \\ 
V^{h}\left( r,\varphi _{i},0\right) =V^{h}\left( r,\varphi _{i},B\right)
=0,i\in \left[ 0,n\right]%
\end{array}%
\right\} .  \label{7.10}
\end{equation}%
In (\ref{7.9}) derivatives $V_{z}^{h},V_{\varphi }^{h}$ are understood as in
(\ref{7.4})-(\ref{7.7}) being supplied by the periodicity condition (\ref%
{7.8}).

\subsection{Completion of the approximate mathematical model}

\label{sec:7.2}

It was stated in the end of section 4 that we construct an approximate
mathematical model for our travel time tomography problem, and then we
verify this model computationally. The first step of this model is the
truncated Fourier-like series (\ref{5.1}). We remind that truncated series
are used quite often in developments of numerical methods for various
inverse problems. On the other hand, convergencies of such methods as the
number of terms tends to infinity are usually not proven, see, e.g. \cite%
{GN,Kab}. This is caused by the ill-posed nature of inverse problems. The
second and the final step of that approximate mathematical model is the
assumption of partial finite differences with condition (\ref{7.1}). This
step is reflected in (\ref{7.2})-(\ref{7.8}).

\section{Convexification}

\label{sec:8}

\subsection{The cost functional of the convexification}

\label{sec:8.1}

Let $M>0$ be an arbitrary number. Let $c>0$ be the number in (\ref{2.1}).
Define the set $B\left( M\right) $ of vector functions $V^{h}\left( r\right) 
$ of (\ref{7.30}) as%
\begin{equation}
\left. 
\begin{array}{c}
B\left( M\right) =\left\{ 
\begin{array}{c}
V^{h}\in H_{N}^{1,h}\left( D^{h}\right) :u^{h}\left( r,\varphi
_{i},z_{j},z_{0}\right) \geq c,z_{0}\in \left[ 0,B\right] \\ 
\left\Vert V^{h}\right\Vert _{H_{N}^{1,h}\left( D^{h}\right) }<M, \\ 
V^{h}\left( r,\varphi _{i},0\right) =\text{ }G_{0}\left( r,\varphi
_{i}\right) ,V^{h}\left( r,\varphi _{i},B\right) =\text{ }G_{B}\left(
r,\varphi _{i}\right) , \\ 
i\in \left[ 0,n\right] ,j\in \left[ 0,k\right] , \\ 
\text{(\ref{7.6})-(\ref{7.8}) hold,}%
\end{array}%
\right\}%
\end{array}%
\right.  \label{8.1}
\end{equation}%
where the space $H_{N}^{1,h}\left( D^{h}\right) $ is defined in (\ref{7.9}).

Thus, we have replaced BVP2 with BVP3:

\textbf{Boundary Value Problem 3 (BVP3). }\emph{Find the vector function }$%
V^{h}\in \overline{B\left( M\right) }.$\emph{\ }

To solve BVP3, we consider the following minimization problem:

\textbf{Minimization Problem}. \emph{Minimize the following functional on
the set }$\overline{B\left( M\right) }:$%
\begin{equation}
J_{\lambda }\left( V^{h}\right) =\left\Vert \left( A_{N}V^{h}+S\left(
V_{\varphi }^{h},V_{z}^{h},r\right) \right) e^{\lambda r}\right\Vert
_{L_{2,N}^{h}\left( D^{h}\right) }^{2}.  \label{8.2}
\end{equation}

Recall that $e^{2\lambda r}$ is the Carleman Weight Function in (\ref{5.140}%
), (\ref{5.14}) for the Volterra integral operator%
\begin{equation*}
\int\nolimits_{r}^{R}f\left( s\right) ds.
\end{equation*}%
Note that we do not use the Tikhonov penalization term \cite{T} in (\ref{8.2}%
), which is rare in the theory of Inverse Problems.

\subsection{The global convergence}

\label{sec:8.2}

We now formulate theorems of the global convergence analysis for the
Minimization Problem. Both their formulations and proofs are similar with
those of \cite{kin3d}. Therefore, we omit proofs. We now briefly explain the
role of the Carleman Weight Function $e^{2\lambda r}$ in the proof of our
central result, which is Theorem 8.1. Given the set $B\left( M\right) $ in (%
\ref{8.1}), this theorem claims the strong convexity of the functional $%
J_{\lambda }\left( V^{h}\right) $ on $\overline{B\left( M\right) }$ for
sufficiently large values of the parameter $\lambda >>1.$ First, consider
for a moment only the quadratic functional $\widetilde{J}_{\lambda }\left(
V^{h}\right) =\left\Vert \left( A_{N}V^{h}\right) e^{\lambda r}\right\Vert
_{L_{2,N}^{h}\left( D^{h}\right) }^{2}.$ Since the matrix $A_{N}$ is
invertible, then this functional is strongly convex on the entire space $%
L_{2,N}^{h}\left( D^{h}\right) $ for any value of $\lambda .$ However, the
presence of the nonlinear term $S\left( V_{\varphi }^{h},V_{z}^{h},r\right) $
in (\ref{8.2}) might destroy the strong convexity property. To dominate this
term, we use in the proof of Theorem 8.1 the Carleman estimate of Theorem
6.1, which, roughly speaking, states that the nonlinear term is dominated by
the strongly convex quadratic term $\left\Vert \left( A_{N}V^{h}\right)
e^{\lambda r}\right\Vert _{L_{2,N}^{h}\left( D^{h}\right) }^{2}.$

Another question is on what exactly \textquotedblleft sufficiently large $%
\lambda "$ is. Philosophically this is very similar with any asymptotic
theory. Indeed, such a theory usually claims that as soon as a certain
parameter $A$ is sufficiently large, a certain formula $B$ is valid with a
good accuracy. In a computational practice, however, one always works with
specific ranges of many parameters. Therefore, only results of numerical
experiments can establish which exactly values of $A$ are sufficient to
ensure a good accuracy of $B$. In particular, we demonstrate below that $%
\lambda =3$ is an optimal value of the parameter $\lambda $ for our
computations. We also note that in all previous publications on the
convexification of this research group optimal values of $\lambda $ were $%
\lambda \in \left[ 1,5\right] $, see, e.g. \cite{KLZ}-\cite{convMFG}.

\textbf{Theorem 8.1 }(strong convexity). \emph{Assume that conditions (\ref%
{1.4})-(\ref{1.8}) hold, and let }$c>0$\emph{\ be the number in (\ref{2.1})
and (\ref{8.1}). Then:}

1. \emph{The functional }$J_{\lambda }\left( V^{h}\right) $\emph{\ has the Fr%
\'{e}chet derivative }$J_{\lambda }^{\prime }\left( V^{h}\right) \in
H_{0,N}^{1,h}\left( D^{h}\right) $\emph{\ at every point }$V^{h}\in 
\overline{B\left( M\right) }$ \emph{and for all }$\lambda >0.$ \emph{\
Furthermore, the Fr\'{e}chet derivative }$J_{\lambda }^{\prime }\left(
V^{h}\right) $ \emph{satisfies Lipschitz continuity condition on }$\overline{%
B\left( M\right) },$\emph{\ i.e. there exists a number }%
\begin{equation*}
\overline{C}=\overline{C}\left( h_{0},c,M,N,D^{h},\lambda \right) >0\emph{\ }
\end{equation*}%
\emph{depending only on listed parameters such that the following estimate
holds:}%
\begin{equation*}
\left\Vert J_{\lambda }^{\prime }\left( V_{2}^{h}\right) -J_{\lambda
}^{\prime }\left( V_{1}^{h}\right) \right\Vert _{L_{2,N}^{h}\left(
D^{h}\right) }\leq \overline{C}\left\Vert V_{2}^{h}-V_{1}^{h}\right\Vert
_{H_{N}^{1,h}\left( \Omega ^{h}\right) },\text{ }\forall
V_{1}^{h},V_{2}^{h}\in \emph{\ }\overline{B\left( M\right) }.
\end{equation*}

2. \emph{There exist a sufficiently large number }$\lambda _{0}\geq 1$ \emph{%
\ and a number} $C>0,$%
\begin{equation}
\lambda _{0}=\lambda _{0}\left( h_{0},c,M,N,D^{h}\right) \geq 1,\text{ }
C=C\left( h_{0},c,M,N,D^{h}\right) >0,  \label{8.3}
\end{equation}%
\emph{both numbers depending} \emph{only on listed parameters, such that for
every }$\lambda \geq \lambda _{0}$\emph{\ the functional }$J_{\lambda
}\left( V^{h}\right) $\emph{\ is strongly convex on the set }$\overline{
B\left( M\right) }.$ \emph{The strong convexity means that the following
estimate holds}$:$\emph{\ \ }%
\begin{equation*}
\begin{split}
J_{\lambda }\left( V_{2}^{h}\right) -J_{\lambda }\left( V_{1}^{h}\right)
-J_{\lambda }^{\prime } &\left( V_{1}^{h}\right) \left(
V_{2}^{h}-V_{1}^{h}\right) \geq C\left\Vert V_{2}^{h}-V_{1}^{h}\right\Vert
_{H_{2,N}^{1,h}\left( \Omega ^{h}\right) }^{2}, \\
& \forall V_{1}^{h},V_{2}^{h}\in \emph{\ }\overline{B\left( M\right) }.
\end{split}%
\end{equation*}

3.\emph{\ Furthermore, for every }$\lambda \geq \lambda _{0}$\emph{\ there
exists unique minimizer }$V_{\min ,\lambda }^{h}\in \ \overline{B\left(
M\right) }$\emph{\ of the functional }$J_{\lambda }\left( V^{h}\right) $%
\emph{\ on the set }$\ \overline{B\left( M\right) }$\emph{\ and the
following inequality holds:}%
\begin{equation*}
J_{\lambda }^{\prime }\left( V_{\min ,\lambda }^{h}\right) \left(
V^{h}-V_{\min ,\lambda }^{h}\right) \geq 0,\text{ }\forall V^{h}\in 
\overline{B\left( M\right) }.
\end{equation*}

Below $C$ denotes different positive numbers depending only on parameters
listed in (\ref{8.3}). In the regularization theory \cite{T}, the minimizer $%
V_{\min ,\lambda }^{h}$ of functional (\ref{8.2}) is called
\textquotedblleft regularized solution". It is important to estimate the
accuracy of the regularized solution depending on the level of the noise in
the data. To do this, we recall first that, following the regularization
theory, we need to assume the existence of the \textquotedblleft ideal"
solution of BVP3, i.e. solution with the noiseless data. The ideal solution
is also called \textquotedblleft exact" solution. We denote this solution $%
V^{h\ast }\in H_{N}^{1,h}\left( D^{h}\right) .$ We denote the noiseless data
in the fourth line of (\ref{8.1}) as $G_{0}^{\ast }\left( r,\varphi
_{i}\right) ,G_{B}^{\ast }\left( r,\varphi _{i}\right) .$ We assume that the
exact solution $V^{h\ast }\in B^{\ast }\left( M\right) ,$ where $B^{\ast
}\left( M\right) $ is the following analog of the set $B\left( M\right) $ in
(\ref{8.1}) 
\begin{equation}
\left. 
\begin{array}{c}
B^{\ast }\left( M\right) = \\ 
=\left\{ 
\begin{array}{c}
V^{h}\in H_{N}^{1,h}\left( D^{h}\right) :u^{h}\left( r,\varphi
_{i},z_{j},z_{0}\right) \geq c,z_{0}\in \left[ 0,B\right] \\ 
\left\Vert V^{h}\right\Vert _{H_{N}^{1,h}\left( D^{h}\right) }<M, \\ 
V^{h}\left( r,\varphi _{i},0\right) =\text{ }G_{0}^{\ast }\left( r,\varphi
_{i}\right) ,V^{h}\left( r,\varphi _{i},B\right) =\text{ }G_{B}^{\ast
}\left( r,\varphi _{i}\right) , \\ 
i\in \left[ 0,n\right] ,j\in \left[ 0,k\right] , \\ 
\text{(\ref{7.6})-(\ref{7.8}) hold.}%
\end{array}%
\right\}%
\end{array}%
\right.  \label{8.4}
\end{equation}

Let $\delta \in \left( 0,1\right) $ be the level of the noise in the
boundary data $G_{0}^{\ast },G_{B}^{\ast }$. To simplify the presentation,
we assume that the noise is not introduced in the function $p\left( \varphi
,z,z_{0}\right) $ in (\ref{1.11}), i.e. we assume that $p\left( \varphi
,z,z_{0}\right) =p^{\ast }\left( \varphi ,z,z_{0}\right) ,$ although this
case can also be included. Suppose that there exist such extensions $%
P^{h}\left( \mathbf{x}^{h}\right) \in H_{N}^{1,h}\left( D^{h}\right) $ and $%
P^{h\ast }\left( \mathbf{x}^{h}\right) \in H_{N}^{1,h}\left( D^{h}\right) $
of the pairs of boundary data $\left( G_{0}\left( r,\varphi _{i}\right)
,G_{B}\left( r,\varphi _{i}\right) \right) _{i=0}^{n}$ and $\left(
G_{0}^{\ast }\left( r,\varphi _{i}\right) ,G_{B}^{\ast }\left( r,\varphi
_{i}\right) \right) _{i=0}^{n}$ that for $i\in \left[ 0,n\right] $ and $r\in
\left( \varepsilon ,R\right) $ 
\begin{equation}
\left. 
\begin{array}{c}
P^{h}\left( r,\varphi _{i},0\right) =G_{0}\left( r,\varphi _{i}\right) ,%
\text{ }P^{h}\left( r,\varphi _{i},B\right) =G_{B}\left( r,\varphi
_{i}\right) , \\ 
P^{h\ast }\left( r,\varphi _{i},0\right) =G_{0}^{\ast }\left( r,\varphi
_{i}\right) ,\text{ }P^{h\ast }\left( r,\varphi _{i},B\right) =G_{B}^{\ast
}\left( r,\varphi _{i}\right) , \\ 
\left\Vert P^{h}\right\Vert _{H_{N}^{1,h}\left( D^{h}\right) }<M,\text{ }%
\left\Vert P^{h\ast }\right\Vert _{H_{N}^{1,h}\left( D^{h}\right) }<M, \\ 
\left\Vert P^{h}-P^{h\ast }\right\Vert _{H_{N}^{1,h}\left( D^{h}\right)
}<\delta .%
\end{array}%
\right. \text{ }  \label{8.5}
\end{equation}

\textbf{Theorem 8.2 }(an estimate of the accuracy of the regularized
solution)\textbf{. }\emph{Assume that conditions of Theorem 8.1 hold. Let
the exact solution }$V^{h\ast }\in B^{\ast }\left( M\right) ,$where the set $%
B^{\ast }\left( M\right) $ \emph{is defined in (\ref{8.4}). Assume that
conditions (\ref{8.5}) hold. Furthermore, assume that }$\left\Vert V^{h\ast
}\right\Vert _{H_{N}^{1,h}\left( \Omega ^{h}\right) }<M-\alpha ,$\emph{\
where the number }$\alpha \in \left( 0,M\right) $\emph{\ is so small that }%
\begin{equation}
\alpha <C\delta .  \label{8.6}
\end{equation}%
\emph{Let }$\lambda _{1}$\emph{\ be the number }$\lambda _{0}$\emph{\ of
Theorem 8.1 in the case when the number }$M$\emph{\ in (\ref{8.3}) is
replaced with }$2M,$%
\begin{equation}
\lambda _{1}=\lambda _{0}\left( h_{0},c,2M,N,D^{h}\right) \geq \lambda
_{0}\left( h_{0},c,M,N,D^{h}\right) .  \label{8.7}
\end{equation}%
\emph{For any }$\lambda \geq \lambda _{1}$ \emph{let }$V_{\min ,\lambda
}^{h}\in \overline{B\left( M\right) }$\emph{\ be the minimizer on of the
functional }$J_{\lambda }\left( V^{h}\right) $\emph{\ on the set }$\overline{%
B\left( M\right) },$\emph{\ which is found in Theorem 8.1. Then the
following accuracy estimate holds:}%
\begin{equation*}
\left\Vert V_{\min ,\lambda }^{h}-V^{h\ast }\right\Vert _{H_{N}^{1,h}\left(
D^{h}\right) }\leq C\delta .
\end{equation*}

For $\lambda \geq \lambda _{1}$ we now construct the gradient descent method
of the minimization of the functional $J_{\lambda }\left( V^{h}\right) .$
Let $\gamma \in \left( 0,1\right) $ be a number and let $V_{0}^{h}\in
B\left( M/3\right) $ be an arbitrary point. The gradient descent method is
constructed as the following sequence:%
\begin{equation}
V_{n}^{h}=V_{n-1}^{h}-\gamma J_{\lambda }^{\prime }\left( V_{n-1}^{h}\right)
,\text{ }n=1,2,...  \label{8.8}
\end{equation}%
Since by Theorem 8.1 $J_{\lambda _{1}}^{\prime }\left( W_{n-1}^{h}\right)
\in H_{0,N}^{1,h}\left( \overline{\Omega }^{h}\right) ,$ $\forall n$, then
all vector functions $V_{n}^{h}$ satisfy the same boundary condition as the
ones in the third line of (\ref{8.1}), see (\ref{7.10}).

\textbf{Theorem 8.3. }\emph{Let }$C\delta \in \left( 0,M/3\right) $\emph{\
in (\ref{8.6}) and let the number }$\beta \in \left( C\delta ,M/3\right) .$%
\emph{\ Suppose that the exact solution }$V^{h\ast }\in B^{\ast }\left(
M/3-\beta \right) $\emph{. Let }$\lambda =\lambda _{1}$\emph{\ where }$%
\lambda _{1}$ \emph{is defined in (\ref{8.7}). Then there exists a
sufficiently small number }$\gamma _{0}\in \left( 0,1\right) $\emph{\ such
that for any }$\gamma \in \left( 0,\gamma _{0}\right) $\emph{\ the sequence
( \ref{8.8}) }$\left\{ V_{n}^{h}\right\} _{n=0}^{\infty }\subset B\left(
M\right) .$\emph{\ In addition, there exists a number }$\rho =\rho \left(
\gamma \right) \in \left( 0,1\right) $\emph{\ such that the following
convergence estimate holds}%
\begin{equation*}
\left\Vert V_{n}^{h}-V^{h\ast }\right\Vert _{H_{N}^{1,h}\left( D^{h}\right)
}\leq C\delta +\rho ^{n}\left\Vert V_{0}^{h}-V_{\min ,\lambda
}^{h}\right\Vert _{H_{N}^{1,h}\left( D^{h}\right) }.
\end{equation*}%
\emph{Furthermore, let the unction }$\left[ n^{h\ast }\left( \mathbf{x}%
^{h}\right) \right] ^{2}$\emph{\ be the exact semi-discrete target function
and\ let }$\left[ n_{n}^{h}\left( \mathbf{x}^{h}\right) \right] ^{2}$\emph{\
be the semi-discrete target function, which is found via the substitution of 
}$V_{n}^{h}$\emph{\ , sequentially, first in the semi-discrete analog of the
first line of (\ref{5.2}), then in (\ref{5.1}) and finally in the
semi-discrete analog of the left hand side of (\ref{3.8}). Then the
following convergence rate holds: }%
\begin{equation*}
\left\Vert \left[ n_{n}^{h}\left( \mathbf{x}^{h}\right) \right] ^{2}-\left[
n^{h\ast }\left( \mathbf{x}^{h}\right) \right] ^{2}\right\Vert
_{L_{2,N}^{h}\left( D^{h}\right) }\leq C\delta +\rho ^{n}\left\Vert
V_{0}^{h}-V_{\min ,\lambda }^{h}\right\Vert _{H_{N}^{1,h}\left( D^{h}\right)
}.
\end{equation*}

Since a smallness assumption is not imposed on the number $M$\ and since the
starting point $V_{0}^{h}$\ of the sequence (\ref{8.8}) is an arbitrary
point of $B\left( M/3\right) ,$\ then Definition of section 1 implies that
Theorem 8.3 ensures the global convergence of the gradient descent method (%
\ref{8.8}).

\section{Numerical Studies}

\label{sec:9}

In this section we describe our numerical studies. We specify parameters of
the domains $\Omega $ and $\Omega _{\varepsilon }$ in (\ref{1.3}) and (\ref%
{1.31}) as: 
\begin{equation*}
R=1,\ \varepsilon =0.01,\ B=1.
\end{equation*}

The refractive index $n(\mathbf{x})$ in eikonal equation (\ref{1.2}) is
taken as 
\begin{equation}
n(\mathbf{x})=\left\{ 
\begin{array}{cc}
c_{a}=const.>1, & \text{inside of the tested inclusion,} \\ 
1, & \text{outside of the tested inclusion.}%
\end{array}%
\right.  \label{9.02}
\end{equation}%
To ensure that $n(\mathbf{x})\in C^{1}\left( \overline{\Omega }\right) $ as
in (\ref{1.7}), we smooth out $n(\mathbf{x})$ near the boundaries of our
tested inclusions. Then we set: 
\begin{align}
& \hspace{1cm}\mbox{correct inclusion/background contrast}=\frac{c_{a}}{1},
\label{9.03} \\
& \text{computed inclusion/background contrast}=\frac{\max_{\text{inclusion}%
}\left( n_{\text{comp}}(\mathbf{x})\right) }{1}.  \label{9.04}
\end{align}%
In the numerical tests below, we take $c_{a}=1.5,3,5$ which means 1.5:1, 3:1
and 5:1 inclusion background contrasts respectively, see (\ref{9.03}). To
demonstrate that our numerical method can work with inclusions, which have
sophisticated non-convex shapes with voids in them, we take shapes
inclusions like letters `$B$' and `$O$', This is similar with our previous
works \cite{KLZ}-\cite{convMFG} on the convexification.

After choosing the function $n(\mathbf{x})$, we use the fast marching
toolbox "Toolbox Fast Marching" \cite{Peyre} in MATLAB to solve eikonal
equation (\ref{1.2}). Then we convert data from Cartesian coordinates to
polar coordinates to gain the observation data in (\ref{1.11})-(\ref{1.13}).
Then we solve the \textbf{Minimization Problem} (\ref{8.2}) to gain the
computed solution $n_{\text{comp}}(r,\varphi ,z)$. Finally, we convert the
reconstructed solution $n_{\text{comp}}(r,\varphi ,z)$ from polar
coordinates to Cartesian coordinates to exhibit results.

To solve eikonal equation (\ref{1.2}) for generating the observation data in
(\ref{1.11})-(\ref{1.13}), we choose $h_{z}=h_{\varphi }=1/40$ in (\ref{7.2}%
) and (\ref{7.3}), as well as we choose $h_{r}=1/40$ to generate discrete
points along the $r-$direction. Because $r\in \left( \varepsilon ,R\right) $
rather than $r\in \left( 0,R\right) $ in (\ref{1.50}), then the first
interval along $r-$direction is $h_{r}-\varepsilon $. To solve the \textbf{%
Minimization Problem} (\ref{8.2}), we choose $h_{z}=h_{\varphi }=h_{r}=1/20$.

To guarantee that the solution of the \textbf{Minimization Problem} (\ref%
{8.2}) satisfies the boundary conditions in (\ref{7.6})-(\ref{7.8}), we
adopt the Matlab's built-in optimization toolbox \textbf{fmincon}. The
iterations of \textbf{fmincon} stop when we get 
\begin{equation}
|\nabla J_{\lambda }\left( V^{h}\right) |<10^{-2}.  \label{9.05}
\end{equation}%
The starting point of iterations of \textbf{fmincon} is chosen as $V^{h}=0$.
Although the starting point does not satisfy the boundary conditions in (\ref%
{7.6})-(\ref{7.8}), \textbf{fmincon }makes sure that the boundary conditions
(\ref{7.6})-(\ref{7.8}) are satisfied on all other iterations of \textbf{%
fmincon}.

We consider the random noise in observation data in (\ref{1.11})-(\ref{1.13}%
) as follows: 
\begin{align}
p^{\xi }\left( \varphi ,z,z_{0}\right) & =p\left( \varphi ,z,z_{0}\right)
\left( 1+\delta \xi _{p}\left( \varphi ,z,z_{0}\right) \right) ,
\label{9.06} \\
p_{0}^{\xi }\left( r,\varphi ,z_{0}\right) & =p_{0}\left( r,\varphi
,z_{0}\right) \left( 1+\delta \xi _{0}\left( r,\varphi ,z_{0}\right) \right)
,  \label{9.07} \\
p_{B}^{\xi }\left( r,\varphi ,z_{0}\right) & =p_{B}\left( r,\varphi
,z_{0}\right) \left( 1+\delta \xi _{B}\left( r,\varphi ,z_{0}\right) \right)
.  \label{9.08}
\end{align}
In (\ref{9.06})-(\ref{9.08}), $\xi _{p}\left( \varphi ,z,z_{0}\right) $ is
the Gaussian random variable depending on variables $\varphi ,z,z_{0}$ and $%
\xi _{0}\left( r,\varphi ,z_{0}\right) ,\xi _{B}\left( r,\varphi
,z_{0}\right) $ are the uniformly distributed random variables in the
interval $[-1,1]$ depending on variables $r,\varphi ,z_{0}$. Also, $\delta
=0.01$ and $\delta =0.03$, which correspond to the $1\%$ and $3\%$ random
noise levels respectively. Recall that, to simplify the presentation, we did
not consider the noise in the function $p\left( \varphi ,z,z_{0}\right) $ in
our theoretical derivations, see (\ref{8.5}). Nevertheless, we still
introduce random noise in this function in our numerical studies, see (\ref%
{9.06}).

To calculate the $\varphi -$derivative and the $z-$derivative of the noisy
function $p^{\xi }( \varphi$, $z$, $z_{0}) $ in (\ref{3.6})-(\ref{3.7}), as
well as the $r-$derivative of noisy functions $p_{0}^{\xi }( r,\varphi
,z_{0})$, $p_{B}^{\xi }( r,\varphi ,z_{0}) $ in (\ref{3.10}), we firstly use
the natural cubic splines to approximate the noisy data in (\ref{9.06})-(\ref%
{9.08}). Then we use the derivatives of those splines to approximate the
derivatives of corresponding noisy data.

The parameters $\lambda $ and $N$ are the two key parameters in our
numerical method. We firstly choose the optimal $\lambda $, when $N$ is
large enough with $N=8$. Then, keeping that optimal $\lambda ,$ we find the
optimal value of $N$. This is done in Test 1. An important point to make
here that, once chosen in Test 1, that optimal pair $\left( \lambda
,N\right) $ is kept for all other tests.

\begin{figure}[tbph]
\centering
\includegraphics[width = 4.5in]{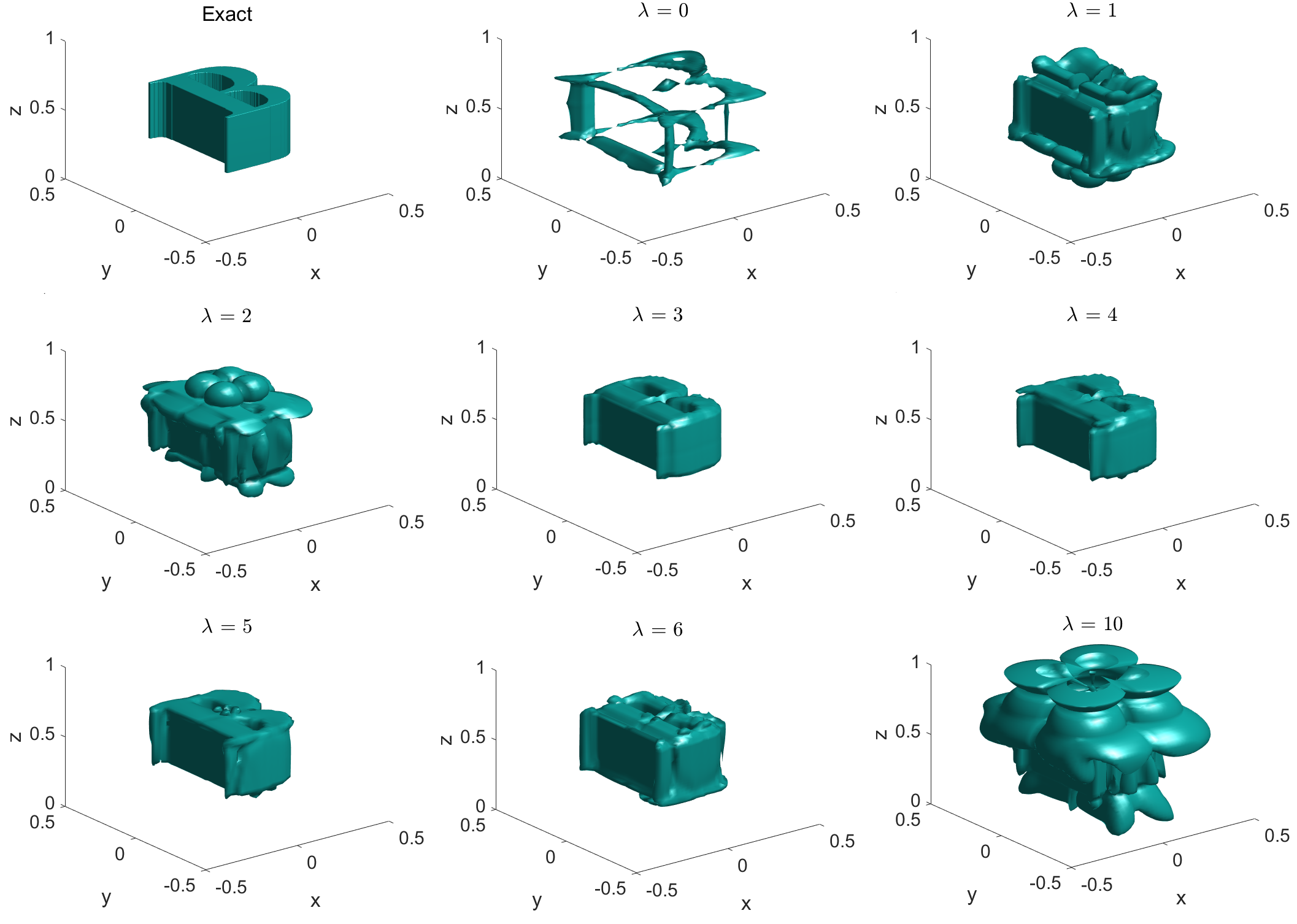}
\caption{Test 1: The reconstructions of $n( \mathbf{x} )$ of different $%
\protect\lambda$ with $N=8$, when the shape of the inclusion in (\protect\ref%
{9.02}) is horizontal letter `$B$' with $c_{a}=1.5$. The value of $\protect%
\lambda $ can be seen on the top side of each square. The images have a low
quality for $\protect\lambda =0,1,2$. The quality is improved with $\protect%
\lambda =3,4,5$. Then it starts to deteriorate at $\protect\lambda =6$ and
becomes unsatisfactory at $\protect\lambda =10$. Thus, we select $\protect%
\lambda=3$ as the optimal value of the parameter $\protect\lambda$. }
\label{diff_lambda}
\end{figure}

\textbf{Test 1.} We test the case when the inclusion in (\ref{9.02}) has the
shape of the horizontally oriented letter `$B$' with $c_{a}=1.5$ in it. The
goal of this test is to find the optimal values of the parameters $\lambda $
and $N$. To make sure that the chosen parameter $N$ does not impact our
results with different values of $\lambda $, we set $N$ to be large enough,
i.e. $N=8$. Computational results are displayed in Figure \ref{diff_lambda}.
We observe that the images have a low quality for $\lambda =0,1,2$. Then the
quality is improved with $\lambda =3,4$, and the reconstruction quality
significantly deteriorates at $\lambda =10$. Hence, we select $\lambda =3$
as the optimal one.

\begin{figure}[tbph]
\centering
\includegraphics[width = 4.5in]{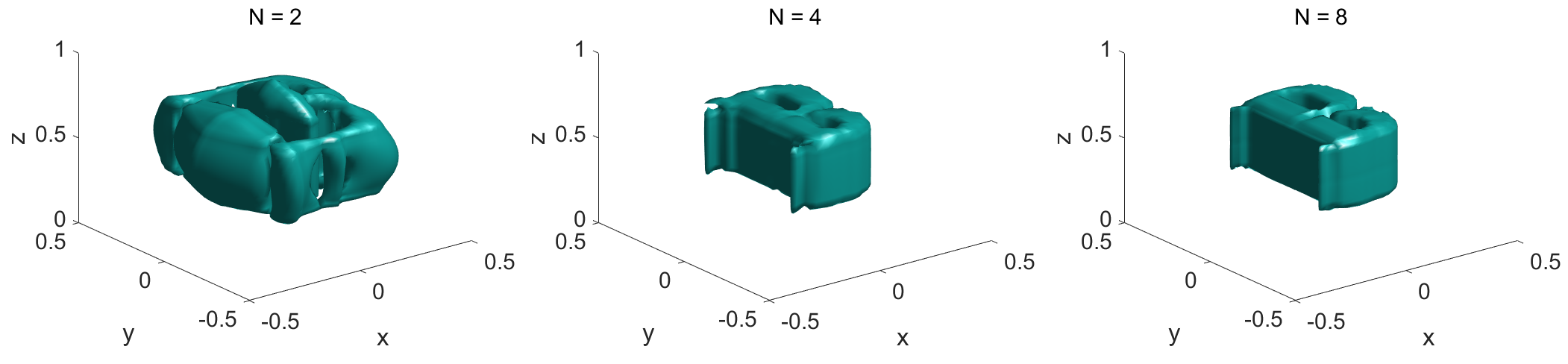}
\caption{Test 1: Reconstructions of $n( \mathbf{x} )$ of $N=2, 4, 8$ at $%
\protect\lambda=3$, when the shape of the inclusion in (\protect\ref{9.02})
is horizontal letter `$B$' with $c_{a}=1.5$. The reconstruction for $N=2$
has a low quality, and the reconstructions for $N=4, 8$ are almost same.
Thus, keeping also in mind (\protect\ref{9.09}) and (\protect\ref{900}) and
also to reduce the computational cost, we choose in Tests 2-5 $N=4$ and $%
\protect\lambda =3$. }
\label{diff_N}
\end{figure}

Next, we fix the optimal value of $\lambda =3$ and test the influence of the
parameter $N$. The results with $N=2,4,8$ are shown in Figure \ref{diff_N}.
We can find that the reconstruction for $N=2$ is of a low quality, and the
reconstructions for $N=4,8$ are almost same. Hence, to reduce the
computational cost we select $N=4$ as the optimal value of $N$, instead of $%
N=8$. Furthermore, we also consider the truncated expansion of the function $%
u\left( r,\varphi ,z,z_{0}\right) $ in (\ref{3.1}) with $N=4$, which is
denoted as $u_{4}\left( r,\varphi ,z,z_{0}\right) $. We have obtained that 
\begin{equation}
\frac{\left\Vert u_{4}\left( r,\varphi ,z,z_{0}\right) \right\Vert
_{L_{2}\left( \Omega \times \left( 0,B\right) \right) }}{\left\Vert u\left(
r,\varphi ,z,z_{0}\right) \right\Vert _{L_{2}\left( \Omega \times \left(
0,B\right) \right) }}=99.86\%.  \label{9.09}
\end{equation}%
It is clear from (\ref{9.09}) that $N=4$ is a quite informative case.

\textbf{Conclusion:} We choose 
\begin{equation}
\lambda =3,N=4  \label{900}
\end{equation}%
as the optimal values of these parameters, and we use these values in Tests
2-5.

\begin{figure}[tbph]
\centering
\includegraphics[width = 4.5in]{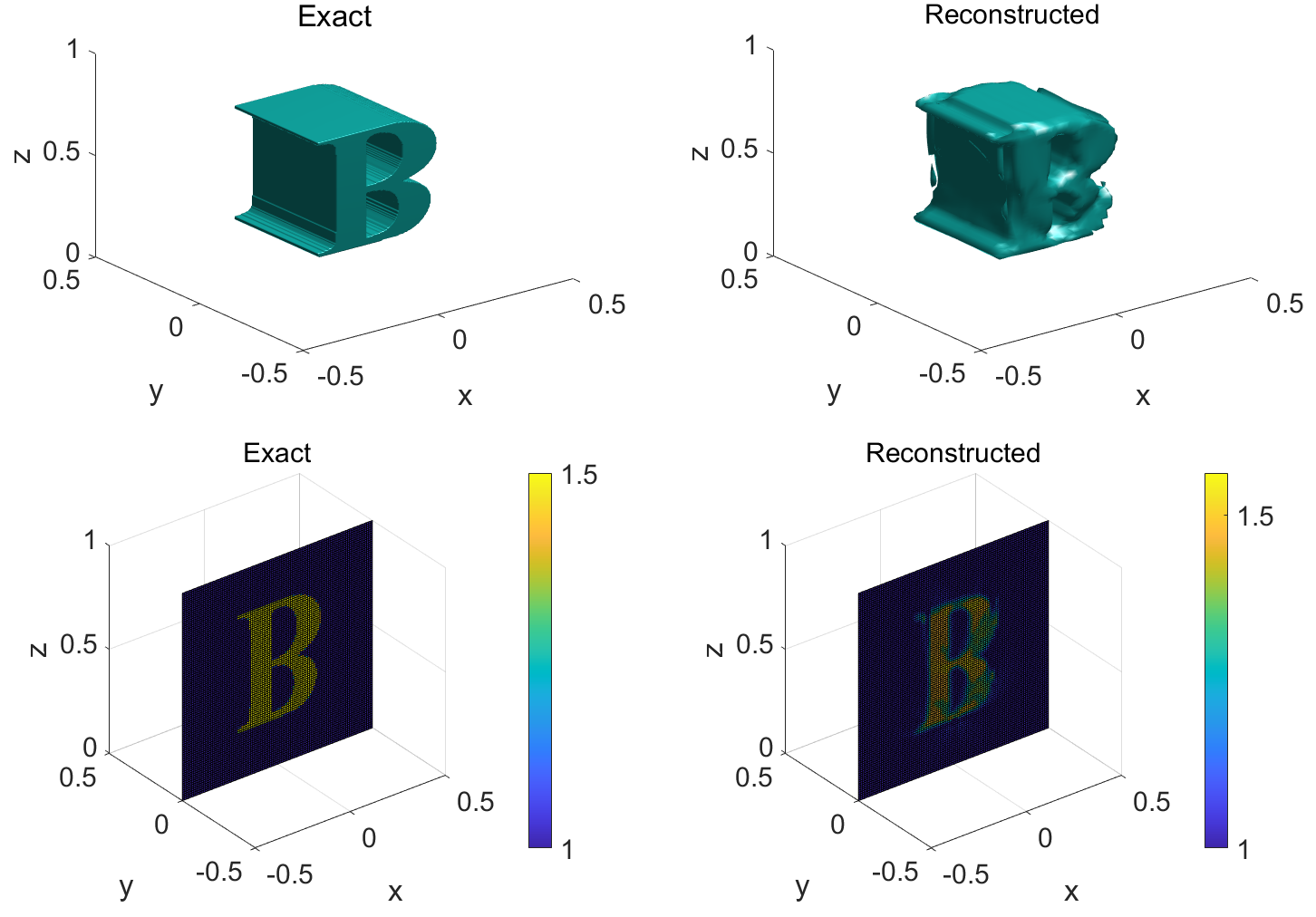}
\caption{Test 2: The exact (left) and reconstructed (right) function $n( 
\mathbf{x}) $, when the shape of the inclusion in (\protect\ref{9.02}) is
vertically oriented letter `$B$' with $c_{a}=1.5$ in it. Here $\protect%
\lambda=3,N=4$ as in (\protect\ref{900}). The reconstruction is accurate. }
\label{plot_re_B02}
\end{figure}

\textbf{Test 2.} We test the case when the inclusion in (\ref{9.02}) has the
shape of the vertically oriented letter `$B$' with $c_{a}=1.5$ in it.
Results are presented on Figure \ref{plot_re_B02}. An accurate
reconstruction can be observed.

\begin{figure}[tbph]
\centering
\includegraphics[width = 4.5in]{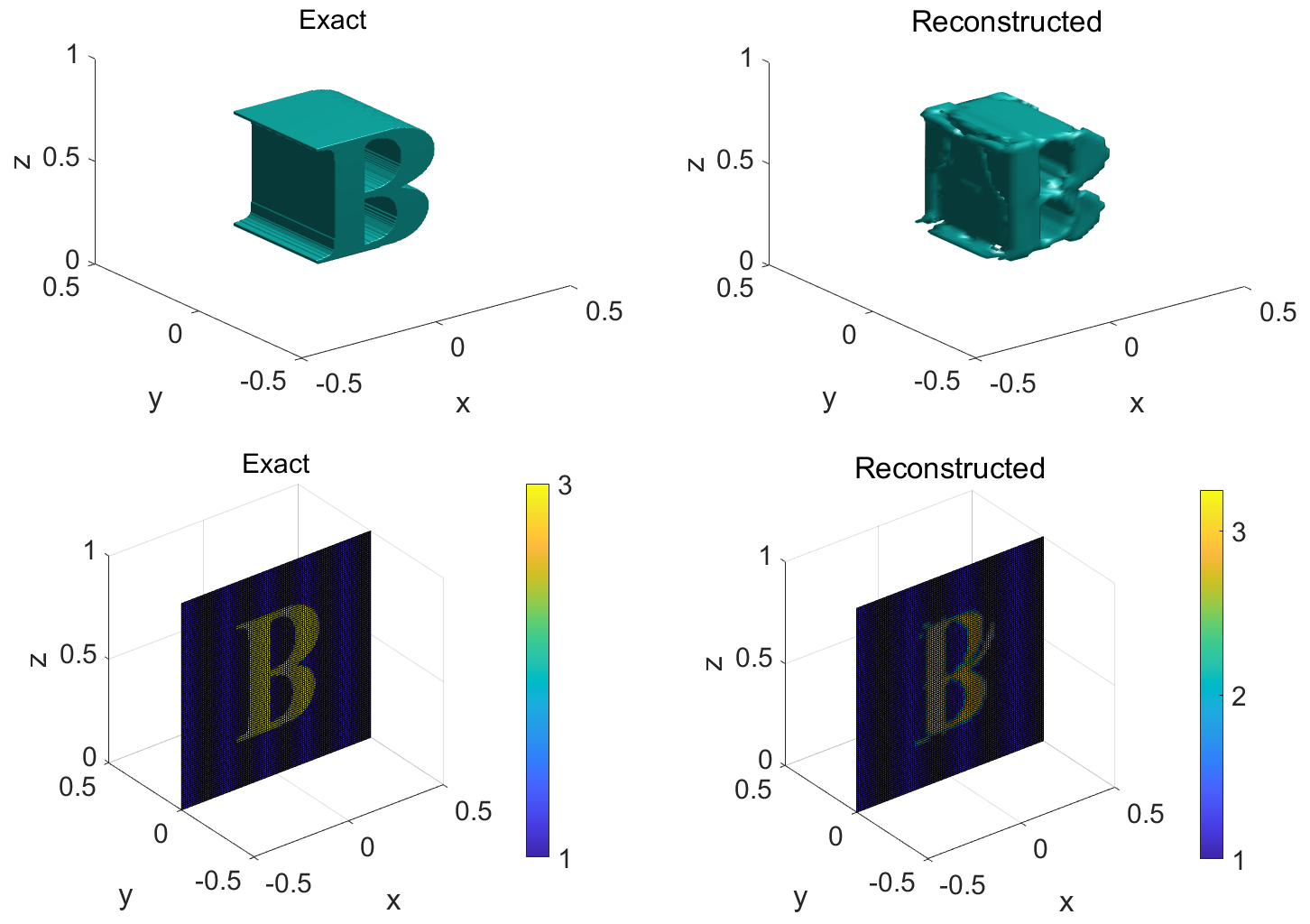}
\caption{Test 3: The exact (left) and reconstructed (right) function $n( 
\mathbf{x}) $, when the shape of the inclusion in (\protect\ref{9.02}) is
vertically oriented letter `$B$' with $c_{a}=3$ in it. Here $\protect\lambda %
=3,N=4$ as in (\protect\ref{900}). The inclusion/background contrast in (%
\protect\ref{9.03}) is $3:1$. The computed inclusion/background contrast in (%
\protect\ref{9.04}) is accurate. }
\label{plot_re_B03}
\end{figure}

\textbf{Test 3.} We test the cases when the inclusion in (\ref{9.02}) has
the shape of the vertically oriented letter `$B$' with $c_{a}=3$ and $%
c_{a}=5 $ in it. Results are displayed on Figures \ref{plot_re_B03} and \ref%
{plot_re_B04} respectively. The inclusion/background contrasts in (\ref{9.03}%
) are respectively $3:1$ and $5:1$. We see that shape of the inclusion is
imaged accurately in both cases. In addition, the computed
inclusion/background contrasts in (\ref{9.04}) are accurate.

\begin{figure}[tbph]
\centering
\includegraphics[width = 4.5in]{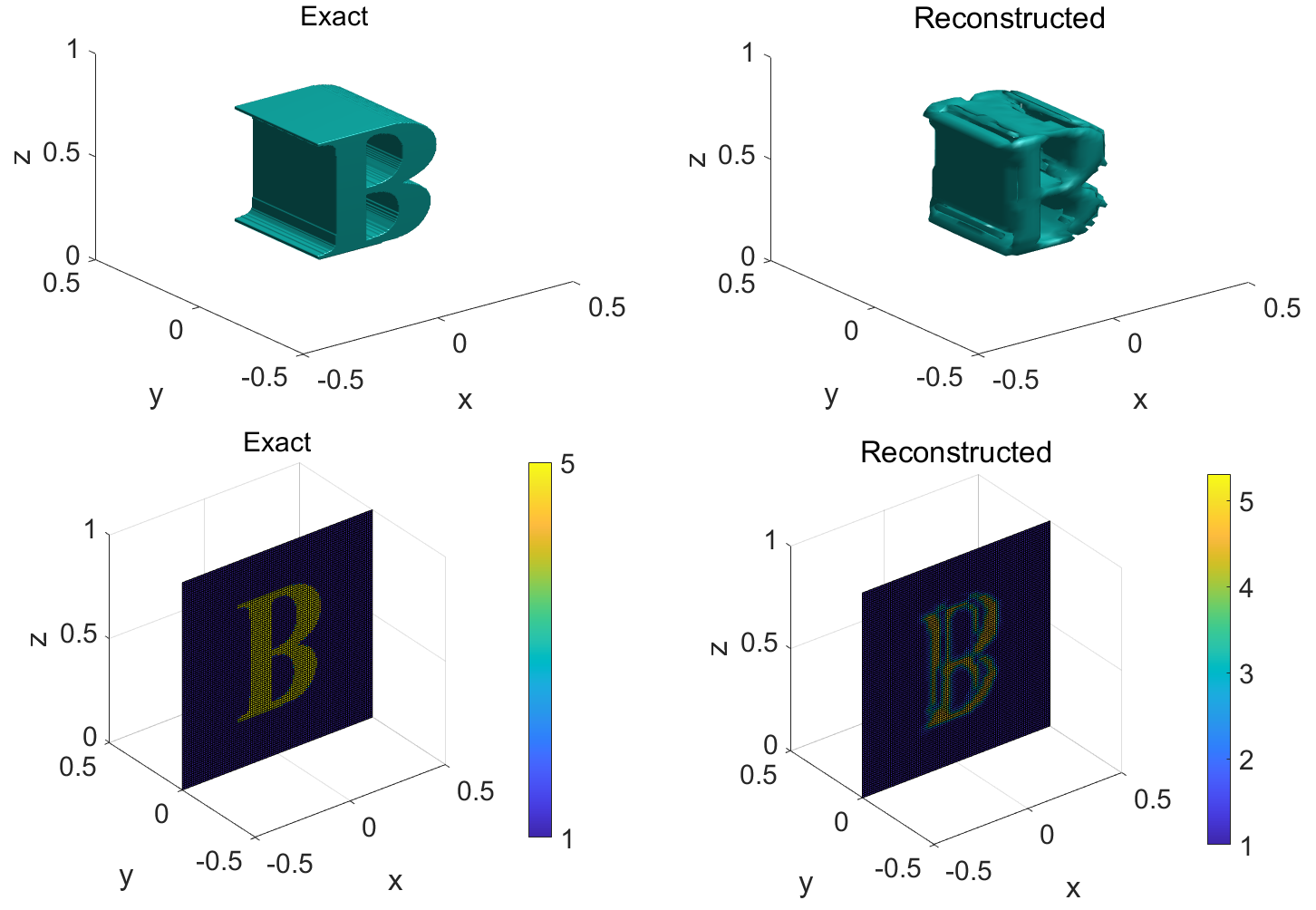}
\caption{Test 3: The exact (left) and reconstructed (right) function $n( 
\mathbf{x}) $, when the shape of the inclusion in (\protect\ref{9.02}) is
vertically oriented letter `$B$' with $c_{a}=5$ in it. Here $\protect\lambda%
=3,N=4$ as in (\protect\ref{900}). The inclusion/background contrast in ( 
\protect\ref{9.03}) is $5:1$. The computed inclusion/background contrast in
( \protect\ref{9.04}) is accurate. }
\label{plot_re_B04}
\end{figure}

\begin{figure}[tbph]
\centering
\includegraphics[width = 4.5in]{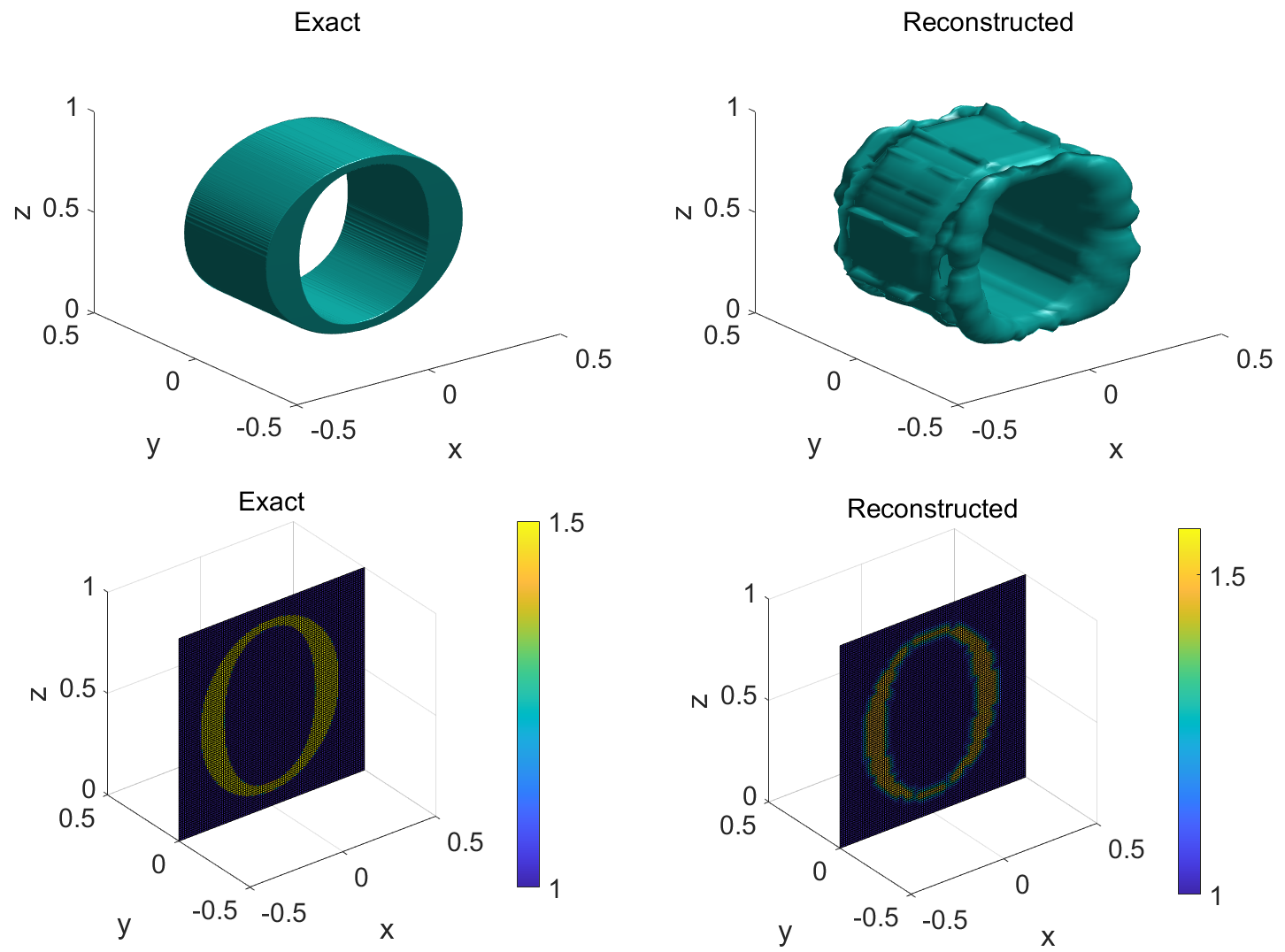}
\caption{Test 4: The exact (left) and reconstructed (right) function $n( 
\mathbf{x}) $, when the shape of the inclusion in (\protect\ref{9.02}) is
the letter `$O$' elongated along the y-axis and with $c_{a}=1.5$ in it. Here 
$\protect\lambda =3,N=4$ as in (\protect\ref{900}). The reconstruction is
accurate. }
\label{plot_re_O}
\end{figure}

\textbf{Test 4.} We test the case when the inclusion in (\ref{9.02}) has the
shape of the letter `$O$' elongated along the $y-$axis and with $c_{a}=1.5$
in it. Results are presented on Figure \ref{plot_re_O}. We again observe an
accurate reconstruction of both the shape of the inclusion and the
inclusion/background contrast.

\begin{figure}[tbph]
\centering
\includegraphics[width = 4.5in]{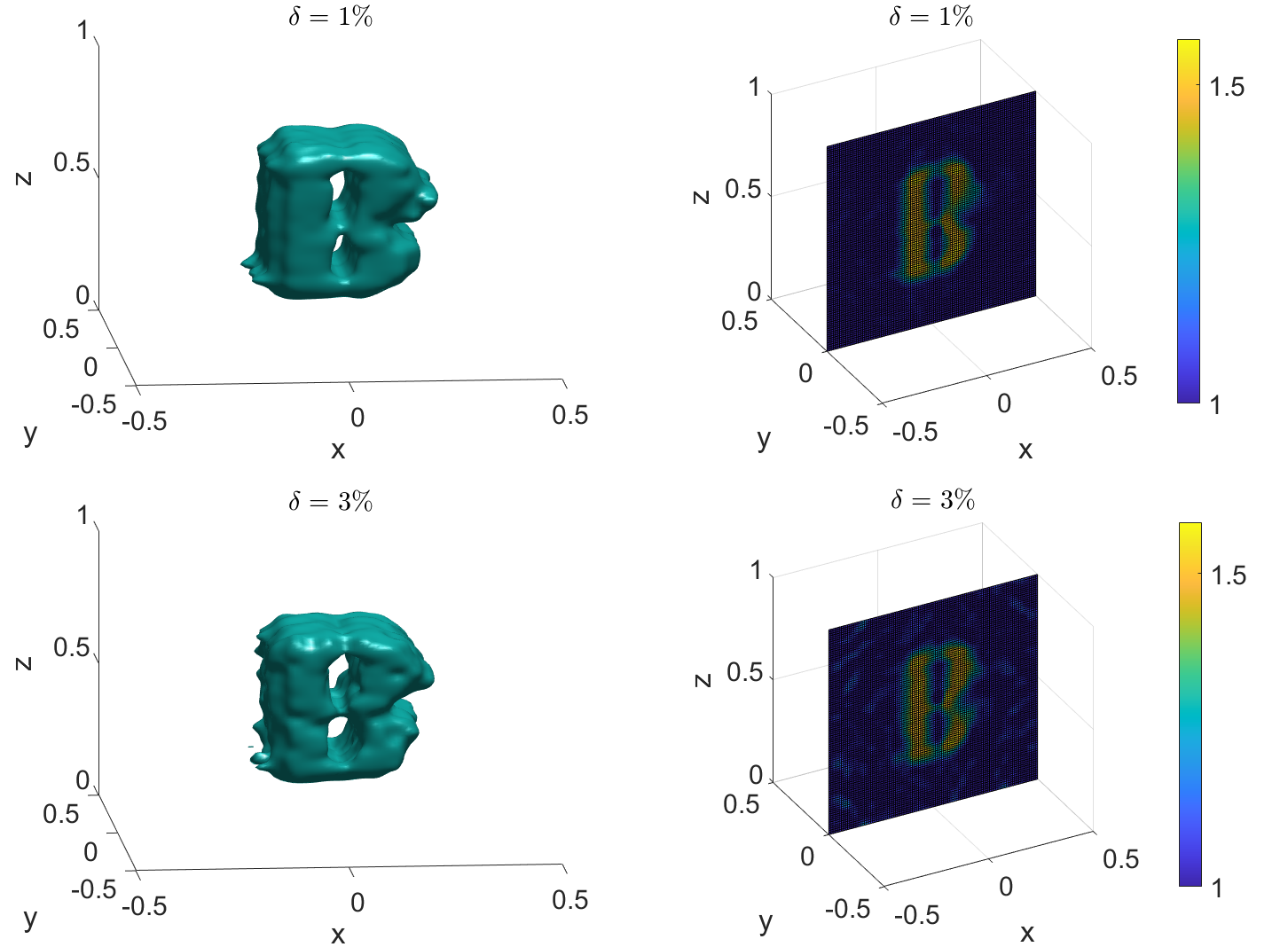}
\caption{Test 5: Reconstructed function $n(\mathbf{x})$ with $\protect\sigma %
=0.01$ (top) and $\protect\sigma =0.03$ (bottom), i.e. with 1\% and 3\%
noise level, when the shape of the inclusion in (\protect\ref{9.02}) is
vertically oriented letter `$B$' with $c_{a}=1.5$ in it. Here $\protect%
\lambda =3,N=4$ as in (\protect\ref{900}). The reconstructions and computed
inclusion/background contrasts in (\protect\ref{9.04}) are accurate. }
\label{plot_noise_B02}
\end{figure}

\textbf{Test 5.} We consider the case when the random noise is present in
the data in (\ref{9.06})-(\ref{9.08}) with $\delta =0.01$ and $\delta =0.03$%
, i.e. with 1\% and 3\% noise level respectively. We test the reconstruction
for the case when the inclusion in (\ref{9.02}) has the shape of the
vertically oriented letter `$B$' with $c_{a}=1.5$ in it. Results are
displayed on Figure \ref{plot_noise_B02}. The reconstructions of the shape
of the inclusion as well as computed inclusion/background contrasts in (\ref%
{9.04}) are accurate.

\textbf{Remarks 9.1:}

\begin{enumerate}
\item \ \emph{Recall that by Theorem 3.1 condition (\ref{1.8}) is a
sufficient condition for our method to work. Recall also that in the data
generation process we smooth out tested inclusions in small neighborhoods of
their boundaries. Given these, a careful analysis of correct images of
Figures 2-8 indicates that condition (\ref{1.8}) is satisfied at least in a
major part of the domain }$\Omega $\emph{\ in each of the above Tests 1-5. }

\item \emph{We conclude, therefore, that our method works numerically under
conditions which are broader than (\ref{1.8}). In other words, (\ref{1.8})
is not a necessary condition for our method to work. Establishing necessary
conditions is outside of the scope of the current publication.}
\end{enumerate}

\textbf{Acknowledgments}. The work of Li was partially supported by
Guangdong Basic and Applied Basic Research Foundation 2023B1515250005. The
work of Romanov was partially supported by the Mathematical Center in
Akademgorodok under Agreement 075152022281 with the Ministry of Science and
Higher Education of the Russian Federation. The work of Yang was partially
supported by Supercomputing Center of Lanzhou University.

%\bibliographystyle{siamplain}
%\bibliography{references-20240913}

\end{document}